
\def\LaTeX{%
  \let\Begin\begin
  \let\End\end
  \let\salta\relax
  \let\finqui\relax
  \let\futuro\relax}

\def\UK{\def\our{our}\let\sz s}
\def\USA{\def\our{or}\let\sz z}

\UK 



\LaTeX

\USA


\salta

\documentclass[twoside,12pt]{article}
\setlength{\textheight}{24cm}
\setlength{\textwidth}{16cm}
\setlength{\oddsidemargin}{2mm}
\setlength{\evensidemargin}{2mm}
\setlength{\topmargin}{-15mm}
\parskip2mm


\usepackage[usenames,dvipsnames]{color}
\usepackage{amsmath}
\usepackage{amsthm}
\usepackage{amssymb,bbm}
\usepackage[mathcal]{euscript}

\usepackage{cite}
%
%
%

\usepackage{hyperref}

\usepackage[ulem=normalem,draft]{changes}
%
%

%
 
\definecolor{viola}{rgb}{0.3,0,0.7}
\definecolor{ciclamino}{rgb}{0.5,0,0.5}
\definecolor{blu}{rgb}{0,0,0.7}
\definecolor{rosso}{rgb}{0.85,0,0}

\def\juerg #1{{\color{viola}#1}}
\def\an #1{{\color{blu}#1}}

\def\last #1{{\color{rosso}#1}}
\def\rev #1{{\color{rosso}#1}}
\def\pier #1{{\color{blue}#1}}

\def\juerg #1{#1}
\def\an #1{#1}
\def\rev #1{#1}
\def\pier #1{#1}
\def\last #1{#1}




\bibliographystyle{plain}


%
\newtheorem{theorem}{Theorem}[section]

\newtheorem{proposition}[theorem]{Proposition}

\finqui

\def\Beq{\Begin{equation}}
\def\Eeq{\End{equation}}

\def\Bthm{\Begin{theorem}}
\def\Ethm{\End{theorem}}
\def\Blem{\Begin{lemma}}
\def\Elem{\End{lemma}}

\def\Bcor{\Begin{corollary}}
\def\Ecor{\End{corollary}}
\def\Brem{\Begin{remark}\rm}
\def\Erem{\End{remark}}

\def\Bdim{\Begin{proof}}
\def\Edim{\End{proof}}
\def\Bcenter{\Begin{center}}
\def\Ecenter{\End{center}}
\let\non\nonumber




\def\step #1 \par{\medskip\noindent{\bf #1.}\quad}
\def\jstep #1: \par {\vspace{2mm}\noindent\underline{\sc #1 :}\par\nobreak\vspace{1mm}\noindent}

\def\aand{\quad\hbox{and}\quad}
\def\Lip{Lip\-schitz}
\def\Holder{H\"older}
\def\Frechet{Fr\'echet}

\def\lhs{left-hand side}
\def\rhs{right-hand side}
\def\sfw{straightforward}



\def\nbh{neighb\our hood}


\def\multibold #1{\def\arg{#1}%
  \ifx\arg\pto \let\next\relax
  \else
  \def\next{\expandafter
    \def\csname #1#1\endcsname{{\boldsymbol #1}}%
    \multibold}%
  \fi \next}

\def\pto{.}

\def\multical #1{\def\arg{#1}%
  \ifx\arg\pto \let\next\relax
  \else
  \def\next{\expandafter
    \def\csname cal#1\endcsname{{\cal #1}}%
    \multical}%
  \fi \next}

\def\multigrass #1{\def\arg{#1}%
  \ifx\arg\pto \let\next\relax
  \else
  \def\next{\expandafter
    \def\csname grass#1\endcsname{{\mathbb #1}}%
    \multigrass}%
  \fi \next}


\def\multimathop #1 {\def\arg{#1}%
  \ifx\arg\pto \let\next\relax
  \else
  \def\next{\expandafter
    \def\csname #1\endcsname{\mathop{\rm #1}\nolimits}%
    \multimathop}%
  \fi \next}

\multibold
qweryuiopasdfghjklzxcvbnmQWERTYUIOPASDFGHJKLZXCVBNM.  

\multical
QWERTYUIOPASDFGHJKLZXCVBNM.

\multigrass
QWERTYUIOPASDFGHJKLZXCVBNM.

\multimathop
diag dist div dom mean meas sign supp .


\def\accorpa #1#2{\eqref{#1}--\eqref{#2}}
\def\Accorpa #1#2 #3 {\gdef #1{\eqref{#2}--\eqref{#3}}%
  \wlog{}\wlog{\string #1 -> #2 - #3}\wlog{}}


\def\separa{\noalign{\allowbreak}}

\def\somma #1#2#3{\sum_{#1=#2}^{#3}}

\def\graffe #1{\mathopen\{#1\mathclose\}}
\def\<#1>{\mathopen\langle #1\mathclose\rangle}
\def\norma #1{\mathopen \| #1\mathclose \|}

\def\aeQ{\checkmmode{a.e.\ in~$Q$}}

\def\aet{\checkmmode{a.e.\ in~$(0,T)$}}

\let\tilde\widetilde
\let\hat\widehat
\def\cpto{\,\cdot\,}

\def\iot {\int_0^t}
\def\ioT {\int_0^T}
\def\intQt{\int_{Q_t}}
\def\intQ{\int_Q}
\def\iO{\int_\Omega}
\def\iG{\int_\Gamma}

\def\dt{\partial_t}
\def\dn{\partial_{\nn}}
\def\ddt{\frac d{dt}}

\let\emb\hookrightarrow
\def\cpto{\,\cdot\,}

\def\checkmmode #1{\relax\ifmmode\hbox{#1}\else{#1}\fi}


\let\erre\grassR
\let\enne\grassN




\def\genspazio #1#2#3#4#5{#1^{#2}(#5,#4;#3)}
\def\spazio #1#2#3{\genspazio {#1}{#2}{#3}T0}

\def\L {\spazio L}
\def\H {\spazio H}

\def\C #1#2{C^{#1}([0,T];#2)}


\def\Lx #1{L^{#1}(\Omega)}
\def\Hx #1{H^{#1}(\Omega)}
\def\Wx #1{W^{#1}(\Omega)}

\def\LQ #1{L^{#1}(Q)}

\def\CS #1{C^{#1}(\Sigma)}

\def\Ldue{\Lx 2}
\def\Linfty{\Lx\infty}

\def\Huno{\Hx 1}
\def\Hdue{\Hx 2}


\def\LLx #1{\LL^{#1}(\Omega)}
\def\HHx #1{\HH^{#1}(\Omega)}

\def\LLQ #1{\LL^{#1}(Q)}


\let\eps\varepsilon
\let\badphi\phi
\let\phi\varphi
\def\zz{\boldsymbol\zeta}

\let\TeXchi\chi                         
\newbox\chibox
\setbox0 \hbox{\mathsurround0pt $\TeXchi$}
\setbox\chibox \hbox{\raise\dp0 \box 0 }
\def\chi{\copy\chibox}



\def\soluz{(\vv,\phi,\mu,w)}

\def\ustar{\uu^*}
\def\vstar{\vv^*}
\def\phistar{\phi^*}
\def\mustar{\mu^*}
\def\wstar{w^*}
\def\soluzstar{(\vstar,\phistar,\mustar,\wstar)}
\def\Lambdastar{\boldsymbol\Lambda^*}

\def\Jtilde{\tilde J}

\def\xxi{\boldsymbol\xi}
\def\oomega{\boldsymbol\omega}

\def\vh{\vv^\hh}
\def\phih{\phi^\hh}
\def\muh{\mu^\hh}
\def\wh{w^\hh}
\def\soluzh{(\vh,\phih,\muh,\wh)}

\def\soluzF{(\xxi,\psi,\eta,\omega)}

\def\soluza{(\oomega,p,q,r)}

\def\un{\uu_n}
\def\vn{\vv_n}
\def\phin{\phi_n}
\def\mun{\mu_n}
\def\wn{w_n}
\def\soluzn{(\vn,\phin,\mun,\wn)}

\def\VVz{\VV_{\!0}}
\def\HHz{\HH_0}
\def\0{\boldsymbol 0}

\let\T\grassT
\let\I\grassI

\def\phiz{\phi_0}

\def\Vp{{V^*}}
\def\Wp{W^*}

\def\normaV #1{\norma{#1}_V}
\def\normaW #1{\norma{#1}_W}

\def\lambdamin{\lambda_*}
\def\lambdamax{\lambda^*}

\def\CK{\calC_K}

\def\CS{\calC_S}
\def\CE{\calC_E}
\def\Cdelta{\calC_\delta}
\def\cdelta{c_\delta}

\def\cM{c_M}

\def\cMN{c_{M,N}}

\def\ck{c_\kappa}
\def\ckM{c_{\kappa,M}}

\def\phiQ{\phi_Q}
\def\phiO{\phi_\Omega}
\def\Uad{\calU_{ad}}
\def\UR{\calU_R}
\def\umin #1{\underline u_{#1}}
\def\umax #1{\overline u_{#1}}
\def\ustar{\uu^*}
\def\gQ{g_Q}
\def\gO{g_\Omega}

\def\LLamu{\an{\boldsymbol{\Lambda}_1}}
\def\LLamd{\an{\boldsymbol{\Lambda}_2}}


\usepackage{amsmath}
\DeclareFontFamily{U}{mathc}{}
\DeclareFontShape{U}{mathc}{m}{it}%
{<->s*[1.03] mathc10}{}

\DeclareMathAlphabet{\mathscr}{U}{mathc}{m}{it}

\Begin{document}


%
\title{
\juerg{Optimal \pier{v}elocity \pier{c}ontrol of a \\Brinkman–Cahn–Hilliard \pier{s}ystem\\ with \pier{c}urvature \pier{e}ffects}
}

\author{}
\date{}
\maketitle
\Bcenter
\vskip-1.5cm
{\large\sc Pierluigi Colli$^{(1)}$}\\
{\normalsize e-mail: {\tt pierluigi.colli@unipv.it}}\\
[0.25cm]
{\large\sc Gianni Gilardi$^{(2)}$}\\
{\normalsize e-mail: {\tt gianni.gilardi@unipv.it}}\\
[0.25cm]
{\large\sc Andrea Signori$^{(3)}$}\\
{\normalsize e-mail: {\tt andrea.signori@polimi.it}}\\
[0.25cm]
{\large\sc J\"urgen Sprekels$^{(4)}$}\\
{\normalsize e-mail: {\tt juergen.sprekels@wias-berlin.de}}\\
[.5cm]
$^{(1)}$
{\small Dipartimento di Matematica ``F. Casorati'', Universit\`a di Pavia}\\
{\small and Research Associate at the IMATI -- C.N.R. Pavia}\\
{\small via Ferrata 5, I-27100 Pavia, Italy}\\
[0.25cm]
$^{(2)}$
{\small Dipartimento di Matematica ``F. Casorati'', Universit\`a di Pavia}\\
{\small via Ferrata 5, I-27100 Pavia, Italy}\\
{\small and Istituto Lombardo, Accademia di Scienze e Lettere}\\
{\small via Borgonuovo 25, I-20121 Milano, Italy}\\
[0.25cm]
$^{(3)}$
{\small Dipartimento di Matematica, Politecnico di Milano}\\
{\small via E. Bonardi 9, I-20133 Milano, Italy}
\\
\an{\small Alexander von Humboldt Research Fellow}
\\[.3cm] 
$^{(4)}$
{\small Department of Mathematics}\\
{\small Humboldt-Universit\"at zu Berlin}\\
{\small Unter den Linden 6, D-10099 Berlin, Germany}\\
{\small and Weierstrass Institute for Applied Analysis and Stochastics}\\
{\small Anton-Wilhelm-Amo-Strasse 39, D-10117 Berlin, Germany}\\[10mm]
\Ecenter
\Begin{abstract}
\noindent
We address an optimal control problem governed by a system coupling a Brinkman-type momentum equation for the velocity field with a sixth-order Cahn--Hilliard equation for the phase variable, incorporating curvature effects in the free energy.
The control acts as a distributed velocity control, allowing \juerg{for} the manipulation of the flow field and, consequently, the phase separation dynamics.
  We establish the existence of optimal controls, prove the Fr\'echet differentiability of the control-to-state operator, and derive first-order necessary optimality conditions \juerg{in terms of a variational  inequality involving the adjoint state variables. We also discuss the aspect of sparsity}.  Beyond its analytical novelty, this work provides 
\pier{a} rigorous control framework for 
Brinkman--Cahn--Hilliard systems incorporating a curvature regularization, offering a foundation for applications in microfluidic design and controlled pattern formation.

\vskip3mm
\noindent {\bf Key words:} 
Brinkman–Cahn–Hilliard system, sixth-order Cahn–Hilliard model, curvature effects, optimal control, Fréchet differentiability, adjoint system, sparsity.

\vskip3mm
\noindent 
{\bf AMS (MOS) Subject Classification:} 35Q35, 
\pier{35M33,} 
%
49K20,  
76D07. 
\End{abstract}

\pagestyle{myheadings}
\newcommand\testopari{\sc Colli -- Gilardi -- Signori -- Sprekels}
\newcommand\testodispari{\sc Optimal control of Brinkman--Cahn--Hilliard systems}
\markboth{\testopari}{\testodispari}
%
\section{Introduction}
\label{INTRO}
\setcounter{equation}{0}
Diffuse-interface models have greatly advanced the study of multiphase flows by providing an energetic framework for interfacial phenomena and curvature effects. In particular, higher-order Cahn--Hilliard systems effectively describe phase separation and membrane dynamics in soft matter and biological contexts. In the recent contribution~\cite{CGSS8}, a Brinkman--Cahn--Hilliard model with curvature effects was introduced, coupling a sixth-order phase field equation with a Brinkman-type momentum balance to capture viscous dissipation and curvature-induced stresses in two-phase incompressible flows.

The present paper aims to take a further step in this line of research by investigating the corresponding optimal control problem. The possibility of influencing the order parameter dynamics through a distributed control \juerg{via} the velocity field provides a natural bridge between theoretical analysis and applications, ranging from microfluidic design to the controlled patterning of soft materials.
Accordingly, we consider here an optimal control problem governed by the aforementioned sixth-order Brinkman--Cahn--Hilliard system, which models phase separation processes in incompressible viscous mixtures and incorporates curvature effects within the free energy functional. The corresponding state system is given by
\begin{alignat}{2}
  & - \div \T(\phi,\vv,p)
  + \lambda(\phi) \vv
  = \mu \nabla\phi
  + \uu
  \aand 
  \div\vv = 0
  && \quad \hbox{in $Q$},
  \label{Iprima}
  \\
  & \dt\phi
  + \vv \cdot \nabla\phi
  - \div(m(\phi)\nabla\mu)
  = S(\phi)
  && \quad \hbox{in $Q$},
  \label{Iseconda}
  \\
  & -\eps \Delta w
  + \tfrac 1 \eps f'(\phi) w
  + \nu w 
  = \mu
  && \quad \hbox{in $Q$},
  \label{Iterza}
  \\
  & -\eps \Delta\phi
  + \tfrac 1 \eps  f(\phi) 
  = w
  && \quad \hbox{in $Q$}.
  \label{Iquarta}  
\end{alignat}
Here, $Q$ \pier{stands for} $\Omega\times(0,T)$, where $\Omega\subset\erre^3$ is the spatial bounded domain with smooth boundary~$\Gamma$ and $T>0$ indicates a fixed final time. 
The system \juerg{is then}  complemented with the boundary and initial conditions
\begin{align}
  & \T(\phi,\vv,p) \nn = \0 
  \aand
  \dn\mu = \dn w = \dn\phi = 0
  & \hbox{on $\Sigma$},
  \label{Ibc}
  \\
  & \phi(0) = \phiz
  &  \hbox{in $\Omega$},
  \label{Icauchy}
\end{align}
\Accorpa\Ipbl Iprima Icauchy
where $\nn$ and $\dn$ denote the outer unit normal vector to $\Gamma$ and the associated normal derivative, respectively,
$\Sigma:=\Gamma\times(0,T)$, and $\phiz$ is a given function acting as initial datum.
\juerg{We observe that, in view of \eqref{Iquarta} and \eqref{Ibc}, the unknown $\phi$ also
(at least formally) satisfies  the boundary
condition $\dn\Delta\phi=0$ on $\Sigma$.}

The first equation \eqref{Iprima} is a Brinkman-type momentum balance, where the \pier{viscosity coefficient $\lambda(\phi)$ in~\eqref{Iprima} is given in terms of a positive function~$\lambda$}, and the Brinkman stress tensor $\T(\phi,\vv,p)$ is defined by
\begin{align*}
	\T(\phi,\vv,p)
	& = \eta(\phi) D \vv
	- p \I 
\end{align*}
with phase-dependent positive viscosity $\eta (\phi)$\juerg{; moreover},
$p$ is the pressure,
$\I\in\erre^{3\times3}$ the identity matrix,
and $D\vv$ denotes the symmetric part of the gradient of the   velocity field~$\vv$,
that is,
\Beq
  D\vv := \frac 12 \, \bigl( \nabla\vv + (\nabla\vv)^\top \bigr) \,.
  \label{symgrad}
\Eeq
This formulation interpolates between the Stokes ($\lambda \equiv 0$) and Darcy ($\eta \equiv 0$) regimes, making it well suited for modeling flows in porous or heterogeneous media.
The variable~$\phi$ serves as the order parameter, representing the local phase concentration, and is normalized so that $\phi = \pm 1$ correspond to the pure phases, while $\juerg{\{-1 < \phi < 1\}}$ describes the diffuse interfacial layer of thickness~$\varepsilon$.
Its evolution is governed by a sixth-order Cahn--Hilliard-type equation with a phase-dependent source term~$S(\phi)$ and a positive mobility \juerg{function}~$m(\phi)$  in equation~\eqref{Iseconda}, whereas~$\uu$ on the \rhs\ of~\eqref{Iprima} denotes the distributed control variable.
Furthermore, in the above equations the quantities $\mu$ and $w$ 
represent the first variations of the total free energy~$\calE(\phi)$ 
and of the Ginzburg--Landau energy, respectively, that is,
$\mu = \frac{\delta\calE}{\delta\phi}$ and $w = \frac{\delta\calG}{\delta\phi}$, with
\begin{align}
 \calE(\phi) 
 & := 
 \calF(\phi) + \nu \calG(\phi)
 =
 \frac 12 \iO  \bigl(-\juerg{\varepsilon}\Delta\phi + \juerg{ \mbox{$\frac 1\eps$}} f(\phi)\bigr)^2
 + \nu  \iO  \bigl( \mbox{$\frac{\juerg{\eps}}2$} \, |\nabla\phi|^2 +  \juerg{ \mbox{$\frac 1\eps$}} F(\phi) \bigr) \,.
  \label{defE} 
\end{align}
In \eqref{defE}, $F$ denotes  a smooth double-well potential,  
$f=F'$,
and $\nu$ is a real parameter, \juerg{which is} not necessarily positive.
A typical example \juerg{for}  $F$ is given by the classical quartic potential
\Beq
  F_{\rm reg}(s) := \frac 14 \, (s^2-1)^2\,, \quad s\in\erre\,.
  \label{Freg}
\Eeq

The well-posedness of problem~\Ipbl\ has been established and discussed in~\cite{CGSS8} within a suitable analytic framework \juerg{that} also includes the analysis of the Darcy limit as the viscosity~$\eta(\phi)$ tends to zero.
The purpose of the present contribution is to build on those results and to provide a rigorous analysis of \juerg{a} corresponding optimal control problem.
For the optimal control application, we consider the tracking-type cost functional
\begin{align*} 
  \calJ(\phi,\uu)
  & {} := \frac{b_1}2 \intQ |\phi-\phiQ|^2
  + \frac{b_2}2 \iO |\phi(T) - \phiO|^2
   + \frac{b_3}2 \intQ |\uu|^2
  +  G(\uu) \,,
\end{align*}
\juerg{which is} to be minimized over the set of admissible controls
\Beq
  \Uad := \graffe{
    \uu=(u_1,u_2,u_3) \in (\juerg{L^\infty(Q)})^3 : \ 
    \umin i \leq u_i \leq \umax i \, \ 
    \hbox{for $i=1,2,3$}
  }\,,
  	\non
\Eeq
subject to the state system~\Ipbl.
In the cost functional, the constants $b_i$ are nonnegative,
and $G$ represents an additional regularization or a sparsity-enhancing term,
\juerg{where a  prototypical choice is} given by
\Beq
\label{Gsparse}
  G(\uu) := \kappa \intQ \juerg{\bigl(|u_1|+|u_2|+|u_3|\bigr)}
\Eeq
\juerg{with some} $\kappa>0$.
In the set $\Uad$ above, the bounds $\umin i$, $\umax i$ are given measurable functions on~$Q$ satisfying
$\umin i \leq \umax i$ \aeQ, ensuring that $\Uad$ is nonempty.
\juerg{At this point, we remark that the notion of sparsity, i.e., the possibility that any locally optimal
control may vanish in subregions of positive measure of the space-time cylinder $Q$,
is very important in practical applications, in particular, in the practical numerical solution of
control problems. In connection with partial differential equations, \pier{sparsity} was first \pier{investigated} in \cite{Stadler}. For an overview \pier{of} the existing literature in this field, we refer the reader to the references given in the recent papers~\cite{SpTr} and~\cite{CoSpTr}.}

\pier{A physical interpretation of the model is discussed in~\cite{CGSS8} (see the references therein).} Here, we briefly recall its relation to diffuse-interface formulations of curvature-driven flows and to applications in optimal control.
A key motivation stems from modeling bilayer membranes and soft-matter systems, where curvature plays a central energetic role.
The classical Helfrich model~\cite{Canham1970, Helfrich1973} describes the elastic bending energy of a smooth membrane $\Gamma_0$ by
\begin{align*}
  {\cal E}_{\rm elastic} = \frac{k}{2} \int_{\Gamma_0} (H - H_0)^2 \, \mathrm{d}S,
\end{align*}
where $H$ is the mean curvature, $H_0$ the spontaneous curvature, and $k$ the bending rigidity.  
An alternative diffuse-interface representation replaces the sharp-interface formulation by introducing an order parameter $\phi$ distinguishing interior ($\phi=1$) and exterior ($\phi=-1$) regions, leading to a modified Willmore functional of the form
\begin{align}\label{Willmore}
  {\cal E}_\eps(\phi)
  = \frac{k}{2\eps} \int_{\Omega}
  \Bigl(-\eps \Delta \phi + \frac{1}{\eps}(\phi^2-1)\phi\Bigr)^2 \,,
\end{align}
where $\eps>0$, as above, represents the interfacial thickness.  
It is known that ${\cal E}_\eps(\phi)$ converges to the sharp-interface energy as $\eps \to 0$  
(see~\cite{DuLiuRyhamWang2005, DuLiuRyhamWang2009}).  
Since sharp-interface asymptotics are not considered here, we fix $\eps=1$ from now on for simplicity.

\juerg{The} energy functional defined in~\eqref{defE}, which characterizes our system, \juerg{can}
 be regarded as a higher-order extension of the classical Ginzburg–Landau free energy~${\cal G}(\phi)$, see~\cite{CGSS8} for further details.
For $\nu=0$, ${\cal E}(\phi)$ reduces to the Willmore-type energy~\eqref{Willmore}, recovering the Canham–Helfrich bending description.
When $\nu>0$, it acts as a curvature-penalized regularization of ${\cal G}(\phi)$, while for $\nu<0$ it yields the functionalized Cahn–Hilliard energy~\cite{GS}, relevant for amphiphilic mixtures and membrane models.
These regimes highlight the flexibility of the framework in describing interfacial phenomena ranging from membrane elasticity to pattern formation.

%

\pier{%
For the aforementioned models, a substantial body of analytical and numerical literature has investigated diffuse-interface formulations.  
Rigorous mathematical analyses can be found in~\cite{SW, ChengWangWiseYuan2020, ClimentEzquerraGuillenGonzalez2019, ColliLaurencot2011, ColliLaurencot2012, WuXu2013, WY, BCMS, GS, Helfrich1973, Canham1970, Ioffe, Kato}, while numerical studies have been carried out in~\cite{AlandEgererLowengrubVoigt2014, CampeloHernandezMachado2006, DuLiuRyhamWang2005, DuLiuRyhamWang2009, DuLiuWang2004, TLVW, LowengrubRatzVoigt2009}.}

\pier{Sixth-order Cahn--Hilliard-type equations have been extensively explored, for instance, in oil\last{-}water\last{-}surfactant dynamics~\cite{PZ1, PZ2, SP} and in phase-field crystal models~\cite{GW1, GW2, M1, M2, WW, CG2}.  
Related optimal control and optimization problems have been studied in~\cite{SigW, CGSS-SIAM, ColliGilardiSprekels2018, FRS, PS, RoccaSprekels2015, SprekelsWu2021, CoSpTr, SprTro, SpTr, Stadler}, whereas the coupling with the standard Cahn--Hilliard equation and Brinkman flows has been analyzed in~\cite{E_thesis, EG1, EG_asy, EL, KS2, CKSS, BCG, ContiG, CGSS8}.}

\pier{Furthermore, these contributions include studies on velocity control for hydrodynamically coupled systems involving either classical or nonlocal Cahn--Hilliard equations~\cite{ColliGilardiSprekels2018, FRS, PS, RoccaSprekels2015, SprekelsWu2021}, as well as sixth-order Cahn--Hilliard systems without fluid coupling~\cite{CGSS-SIAM}.  
Overall, this literature highlights the richness of both analytical and numerical approaches in high-order diffuse-interface models and their optimal control applications, encompassing vesicle-fluid interactions, phase-field crystal dynamics, and complex multiphase flows.}

The nonlinear coupling between the flow and the sixth-order Cahn–Hilliard subsystem, together with the curvature term in the energy~\eqref{defE}, makes the analysis of the problem particularly challenging.
The paper is organized as follows.
In Section~\ref{STATEMENT}, we introduce the notation, assumptions, and main theoretical results.
Section~\ref{CONTROL} is devoted to proving the existence of at least one optimal control for the considered minimization problem.
In Section~\ref{FRECHET}, we study the differentiability properties of the control-to-state operator and establish its Fréchet differentiability in a suitable functional framework. As a preliminary step, we analyze the corresponding linearized system.
In Section~\ref{ADJOINT}, we derive the first-order necessary conditions for optimality, which are then reformulated in terms of an adjoint system that we introduce and solve.
\juerg{The final Section~\ref{SPARSITY} then brings the derivation of sparsity results for the prototypical
sparsity  \last{function $G$ given in} \eqref{Gsparse}.}
Overall, these results constitute a first step toward a rigorous optimal control theory for Brinkman–Cahn–Hilliard systems with curvature effects and lay the groundwork for further analytical and numerical studies.

\section{Notation, assumptions and main results}
\label{STATEMENT}
\setcounter{equation}{0}
In this section, we fix our notation, list the assumptions on the state system,
give it a precise form,
and recall some results for it that are already known.
As for the set $\Omega$, its boundary~$\Gamma$, the normal derivative~$\dn$, 
the space-time cylinder~$Q$ and its lateral surface~$\Sigma$, 
we keep the notation used in the Introduction.
More precisely, $\Omega$~is an open, \juerg{bounded and connected} set in $\erre^3$ 
\juerg{having a smooth boundary}.
The symbol $|\Omega|$ denotes its Lebesgue measure.
For any Banach space~$X$, the symbols $\norma\cpto_X$ and $X^*$ denote 
the corresponding norm and its dual space.
However, some exceptions \juerg{to} this notation are listed below.
We introduce the spaces
\begin{align}
  & H := \Ldue , \quad  
  V := \Huno , \quad
  W := \graffe{z\in\Hdue: \ \dn z=0 \ \hbox{on $\Gamma$}} 
  \label{defspazia}
  \\
  & \HH := H \times H \times H \,, \quad
  \VV := V \times V \times V
  \aand
  \VVz := \graffe{\zz\in\VV : \ \div\zz=0}.
 \label{defspazib}
\end{align}
\Accorpa\Defspazi defspazia defspazib
Similarly, we use the boldface characters, like $\LLx2$ and~$\HHx1$, to denote powers of the Lebesgue and Sobolev spaces.
To~simplify the notation, the norms in the special cases $H$ and $\HH$ are simply indicated by~$\norma\cpto$.
The symbol $\norma\cpto_\infty$ might denote the norm in each of the spaces 
$\Linfty$, $\LQ\infty$ and $L^\infty(0,T)$, if no confusion can arise.
Furthermore, we use the same symbol for the norm in some space and the norm in any power 
\juerg{thereof}.

Since $W$ is dense in~$V$, and $V$ is dense in~$H$, 
we can make standard identifications and adopt the usual framework of Hilbert triplets.
Namely, we have that
\begin{align}
  & \< y,z >_V = \textstyle\iO y z 
  \aand
  \< y,z >_W = \< y,z >_V
  \non
  \\
  & \quad \hbox{for every $y\in H$ and $z\in V$ and every $y\in\Vp$ and $z\in W$, respectively}
  \non
\end{align}
so that
\Beq
  W \emb V \emb H \emb \Vp \emb \Wp.
  \label{embeddings}
\Eeq
The symbols $\<\cpto,\cpto>_V$ and $\<\cpto,\cpto>_W$ denote 
the duality pairings between $\Vp$ and~$V$ and between $\Wp$ and~$W$, respectively.
However, thanks to the above identifications, 
there is no confusion if we avoid using subscripts and adopt the simpler symbol $\<\cpto,\cpto>$ for the above pairings.
Besides the space $\VVz$ already introduced, we also define 
\Beq
  \HHz := \graffe{\zz\in\HH : \ \div\zz=0},
  \label{defHHz}
\Eeq
where the divergence is understood in the sense of distributions.
\juerg{Notice}  that the embedding
\Beq
  \VVz \emb \HHz 
  \label{embeddingsbis}
\Eeq
is dense (see, e.g., \cite[Cor.~2.3]{Kato})
and that all of the embeddings in \eqref{embeddings} and \eqref{embeddingsbis} are compact. 
Next, we recall the symbol $D\vv$ introduced in \eqref{symgrad} for the symmetrized gradient of the velocity~$\vv$,
whose use will be extended to any vector field $\zz\in\VV$.
Finally, we adopt the standard notation
\Beq
  A:B := \somma {i,j}13 a_{ij} b_{ij}
  \aand
  |A|^2 := A:A ,
  \quad \hbox{for $A=(a_{ij}),\,B=(b_{ij})\in\erre^{3\times3}$},
  \label{matrices}
\Eeq
for the scalar product and the norm of matrices.

We are ready to list our assumptions on the structure of the state system.
\juerg{For the involved functions and parameters, we postulate:}
\begin{align}
  & \hbox{$\lambda:\erre\to\erre$ is of class $C^1$ and satisfies}
  \quad \lambdamin \leq \lambda(s) \leq \lambdamax \quad
  \non
  \\
  & \quad \hbox{for every $s\in\erre$ and some positive constants $\lambdamin$ and $\lambdamax$}.
  \label{hplambda}
  \\
  & \an{\hbox{$\eta(\cdot)\equiv \eta_0$, $m(\cdot)\equiv m$  with $\eta_0,m>0$, and $\nu$ is real constant}.}
  \label{hpcoeff}
  \\
  & S(s) := -\sigma s + h(s)
  \quad \hbox{for $s\in\erre$, \quad with} \quad
  \sigma \in \erre
  \quad\mbox{and a function}
  \non
  \\
  & \quad \hbox{$h:\erre\to\erre$ of class $C^2$, which is bounded and \Lip\ continuous}.
  \label{hpS}
  \\
  & F : \erre \to \erre
  \quad \hbox{is of class $C^4$, and $F$ and $f:=F'$ satisfy}
  \label{hppot}
  \\
  & \lim_{|s|\to+\infty} \frac{f(s)}s = + \infty\,, 
  \quad \hbox{as well as}
  \label{hpf}
  \\
  & f'(s) \geq -C_1 \,, \quad
  |F(s)| \leq C_2 \, (|s f(s)| + 1)\,,
  \aand
  |s f'(s)| \leq C_3 (|f(s)| + 1)\,,
  \non
  \\
  & \quad \hbox{for every $s\in\erre$ and some positive constants $C_1$, $C_2$ and $C_3$}\,.
  \label{hpFf}
\end{align}
\Accorpa\Hpstruttura hplambda hpFf
Notice that our assumptions on the potential are satisfied if $F=F_{reg}$,
the classical regular potential defined in~\eqref{Freg}.
For the data, we assume that
\begin{align}
  & \uu \in \L2\HH ,
  \label{regu}
  \\
  & \phiz \in W \,.
  \label{regphiz}
\end{align}
\Accorpa\Regdati regu regphiz
As \juerg{stated} in the Introduction, a well-posed\rev{ne}ss result was proved in \cite{CGSS8} 
for an equivalent problem in which the pressure $p$ \juerg{no longer appears}.
Namely, in the quoted paper it was shown that, instead of looking for the pair $(\vv,p)$ solving \eqref{Iprima},
we can equivalently look for its first component $\vv$ as a divergence-free solution 
to the variational equation
\begin{align}
  &\iO \bigl( \eta_0 D\vv : \nabla\zz + \lambda(\phi) \vv \cdot \zz \bigr) 
  = \iO (\mu \nabla\phi + \uu) \cdot \zz\non\\
  &\quad{} \hbox{for every $\zz\in\VVz$ and \aet} \,.
  \label{laxmil}
\end{align}
The argument was based on the identity 
\Beq
  D\zz : \nabla\zz = |D\zz|^2
  \quad \hbox{for every $\zz\in\VV$},
  \label{identity}
\Eeq
and the Korn inequality
\Beq
  \normaV\zz^2 \leq \CK \iO (|D\zz|^2 + |\zz|^2)
  \quad \hbox{for every $\zz\in\VV$} ,
  \label{korn}
\Eeq
which holds true for some constant $\CK>0$ depending only on~$\Omega$.
The combination of these facts yields the coerciveness inequality
\begin{align}
  & \iO \bigl( \eta_0 D\zz : \nabla\zz + \lambda(\badphi) |\zz|^2 \bigr) 
  \geq \alpha \, \normaV\zz^2
  \quad \hbox{for every $\zz\in\VVz$}\,,
  \label{coercive}
  \\
  & \hbox{where} \quad 
  \alpha := \frac{\min\graffe{\eta_0,\lambdamin}}\CK \,, 
  \label{defalpha}
\end{align}
which allows to apply the Lax--Milgram theorem.
Next, \juerg{following the lines of} \cite{CGSS8}, we present another possible version of the problem.
Namely, one can eliminate $w$ by \juerg{inserting \eqref{Iquarta} in} \eqref{Iterza}  \pier{(where now $\eps=1$ in both equations)} to obtain 
\Beq
  - \Delta \bigl( -\Delta\phi + f(\phi) \bigr)
  + \bigl( f'(\phi) + \nu \bigr) \bigl( -\Delta\phi + f(\phi) \bigr)
  = \mu \quad \hbox{in $Q$.}
  \label{quintastrong}
\Eeq

Here is the precise formulation of the problem under consideration.
We look for a quadruplet $\soluz$ with the regularity
\begin{align}
  & \vv \in \L2\VVz\,,
  \label{regv}
  \\
  & \phi \in \H1\Vp \cap \L\infty W \cap \L2{\Hx5} \,,
  \label{regphi}
  \\
  \separa
  & \mu \in \L2V  \,,
  \label{regmu}
  \\
  & w \in \L\infty H \cap \L2{\Hx3\cap W} \,,
  \label{regw}
\end{align}
\Accorpa\Regsoluz regv regw
that solves the variational equations
\begin{align}
  & \iO \bigl(
    \eta_0 D\vv : \nabla\zz
    + \lambda(\phi) \vv \cdot \zz
  \bigr)
  = \iO (\mu \nabla\phi + \uu) \cdot \zz
  \non
  \\
  & \quad \hbox{for every $\zz\in\VVz$ and \aet}\,,
  \label{prima}
  \\
  \separa
  & \< \dt\phi , z > 
  + \iO \vv \cdot \nabla\phi \, z
  + \iO m \nabla\mu \cdot \nabla z
  = \iO S(\phi) z
  \non
  \\
  & \quad \hbox{for every $z\in V$ and \aet}\,,
  \label{seconda}
  \\
  \separa
  &
  \an{\iO \nabla w \cdot \nabla z}
  \pier{{}+{}} \iO \bigl( f'(\phi) + \nu \bigr) w z 
  = \iO \mu z
  \non
  \\
  & \quad \hbox{for every $\an{z\in V}$ and \aet}\,,
  \label{terza}
  \\
  \separa
  & \iO \nabla\phi \cdot \nabla z
  + \iO f(\phi) z
  = \iO w z
  \non
  \\
  & \quad \hbox{for every $z\in V$ and \aet}\,,
  \label{quarta}
\end{align}
as well as the initial condition
\Beq
  \phi(0) = \phiz \,.
  \label{cauchy}
\Eeq
\Accorpa\Pbl prima cauchy
Observe that \eqref{terza} and \eqref{quarta} can be replaced by their strong \juerg{forms}\pier{%
\begin{alignat}{2}
  & -\Delta w +  f'(\phi) w
  + \nu w 
  = \mu
  \quad &&\hbox{ a.e. in $Q$},
  \label{P-terza}
  \\
  & -\Delta\phi
  +  f(\phi) 
  = w
 \quad &&\hbox{ a.e. in $Q$}\juerg{,}
  \label{P-quarta}  
\end{alignat}
thanks} to the regularity of $w$ and $\phi$ required in \eqref{regw} and~\eqref{regphi}
(which also \juerg{encode} the homogeneous Neumann boundary condition for these components).

\pier{Taking \eqref{P-quarta} into account, we can rewrite~\eqref{terza} as}
\begin{align}
  &
  \int_{\Omega} \an{\nabla\!\left( -\Delta \phi + f(\phi) \right) \cdot \nabla z}
  + \int_{\Omega} \bigl(f'(\phi) + \nu\bigr)\,
      \bigl( -\Delta \phi + f(\phi) \bigr) z
  = \int_{\Omega} \mu\, z
  \nonumber
  \\
  & \qquad\text{for every } z \in V \text{ and a.e.\ in } (0,T),
  \label{quinta}
\end{align}
\pier{and keep \eqref{quarta} or \eqref{P-quarta} as the definition of $w$.  
Note that \eqref{quinta} provides a weak formulation of~\eqref{quintastrong}.}

Here are the already known well-posedness and continuous dependence results 
(see \cite[Thm.~2.1 and Thm.~2.3]{CGSS8}).

\Bthm
\label{Wellposedness}
Assume \Hpstruttura\ on the structure, and suppose that the data satisfy \Regdati.
Then there exists a unique quadruplet $\soluz$ with the regularity \Regsoluz\
that solves Problem \Pbl.
Moreover, this solution satisfies the estimate
\begin{align}
  & \norma\vv_{\L2\VV}
  + \norma\phi_{\H1\Vp\cap\L\infty W\cap\L2{\Hx5}}
  \non
  \\
  & \quad {}
  + \norma\mu_{\L2V}
  + \norma w_{\L\infty H\cap\L2{\Hx3}}
  \leq K_1 \,,
  \label{stability}
\end{align}
with a constant $K_1$ that depends only on the structure of the system, $\Omega$, $T$ 
and an upper bound for the norms of the data related to \Regdati.
\Ethm

\Bthm
\label{Contdep}
Assume \Hpstruttura\ for the structure \juerg{and} 
\eqref{regphiz} for the initial datum.
\juerg{Then the following holds true:} if $\an{\uu_i \in \L2\HH}$, $i=1,2$, 
are given and $(\vv_i,\phi_i,\mu_i,w_i)$ are the corresponding solutions,
\juerg{then} the estimate
\begin{align}
  & \norma\vv_{\L2\VV}
  + \norma\phi_{\C0V\cap\L2{\Hx4}}
  + \norma\mu_{\L2H}
  \non
  \\
  & \quad {}
  + \norma w_{\L2W}
  \leq K_2 \, \norma\uu_{\L2\HH}
  \label{contdep}
\end{align}
holds true for the differences $\soluz=(\vv_1,\phi_1,\mu_1,w_1)-(\vv_2,\phi_2,\mu_2,w_2)$ and $\uu=\uu_1-\uu_2$,
with a constant $K_2$ that depends only on the structure of the system, $\Omega$, $T$, the initial datum~$\phiz$,
and an upper bound for the norms of $\uu_1$ and $\uu_2$ in~$\L2\HH$.
\Ethm

Once well-posedness of the state system is established,
one can deal with the control problem.
We present the cost functional mentioned in the Introduction in a precise form.
We~set
\begin{align} 
  & \calJ(\phi,\uu)
  := J(\phi,\uu) + G(\uu)\,,
  \quad \hbox{where}
  \label{cost}
  \\
  & J(\phi,\uu)
  := \frac{b_1}2 \intQ |\phi-\phiQ|^2
  + \frac{b_2}2 \iO |\phi(T) - \phiO|^2
  + \frac{b_3}2 \intQ |\uu|^2\,,
  \label{standard}
  \\
  & \mbox{and}\quad G : \LLQ2 \to [0,+\infty] \quad
  \hbox{is convex, proper and lower semicontinuous}.
  \label{hpG}
\end{align}
In \eqref{standard}, the coefficients and the target functions satisfy
\begin{align}
  & b_1 \,,\, b_2 \in [0,+\infty)
  \aand
  b_3 \in (0,+\infty) ,
  \label{hpbi}
  \\
  & \phiQ \in \LQ2
  \aand
  \phiO \in \last{H} \,.
  \label{hpphiQO}
\end{align}
\Accorpa\Hpcost cost hpphiQO
The set of admissible controls is defined by
\begin{align}
  & \Uad := \graffe{
    \uu=(u_1,u_2,u_3) \in \an{\LL^\infty(Q)} : \ 
    \umin i \leq u_i \leq \umax i \ 
    \hbox{ a.e. in $Q$\,, for $i=1,2,3$}
  },
  \label{Uad}
  \\
  & \quad \hbox{where} \quad
  \umin i \,,\, \umax i \in\LQ\infty
  \quad \hbox{satisfy} 
  \non
  \\
  & \quad
  \umin i \leq \umax i
  \quad \hbox{\aeQ\,, \ for $i=1,2,3$} \,.
  \label{hpuminmax}
\end{align}
\Accorpa\HpUad Uad hpuminmax
\Accorpa\Hpcontrol cost hpuminmax
Then, the control problem \juerg{under investigation} reads:
\begin{align}
  & \hbox{Minimize}\quad \calJ(\phi,\uu)\quad \hbox{over\, $\Uad$}
  \quad \hbox{under the constraint:}
  \non
  \\
  & \quad \hbox{$\phi$ is the second component of the solution}
  \non
  \\
  & \quad \hbox{\pier{to the problem \Pbl\ corresponding to $\uu$}} \,.
  \label{control}
\end{align}
\last{We are ready to present the existence theorem for an optimal control.}

\Bthm
\label{Existencecontrol}
Assume \Hpstruttura\ for the structure of the state system 
and \eqref{regphiz} for the initial dat\an{a}.
Moreover, assume \Hpcost\ for the cost functional and \HpUad\ for the set of the admissible controls.
Then there exists at least \juerg{one} optimal control to the control problem \eqref{control},
that is, there \juerg{is some}  $\ustar\in\Uad$ such that 
$\calJ(\phistar,\ustar)\leq\calJ(\phi,\uu)$ for every $\uu\in\Uad$,
where $\phistar$ and $\phi$ are the second components of the solutions
$\soluzstar$ and $\soluz$ to the state systems \Pbl\ corresponding to $\ustar$ and~$\uu$, respectively.
\Ethm

\rev{The
next step consists} in finding a significant necessary condition for a given admissible control to be optimal.
To this end, we follow a standard procedure.
Namely, we introduce the so-called {\it control-to-state operator}, termed~$\calS$, 
that maps some \nbh\ $\UR$ of $\Uad$ \juerg{in $L^2(0,T;\HH)$} to a suitable Banach space 
and associates to every $\uu\in\UR$ the corresponding state, i.e., 
the solution $\soluz$ to problem \Pbl.
Then, we introduce the composite map $\Jtilde$\an{, also referred to as the {\it reduced cost functional},} by setting, for $\uu\in\UR$,
\Beq
  \Jtilde(\uu) := J(\phi,\uu)\,,
  \quad \hbox{where $\phi$ is the second component of $\soluz:=\calS(\uu)$}\,, 
  \label{defJtilde}
\Eeq
so that the functional to be minimized on $\Uad$ is just $\Jtilde+G$.
Then, a standard argument from Convex Analysis ensures that, 
if $\ustar$ is an optimal control, then there exists some
$\Lambdastar$ in the subdifferential $\partial G(\ustar)\subset\LLQ2$ such that
\Beq
  D\Jtilde(\ustar)[\uu-\ustar] 
  + \intQ \Lambdastar \cdot (\uu-\ustar)
  \geq 0
  \quad \hbox{for every $\uu\in\Uad$}\,,
  \label{preNC}
\Eeq
where $D\Jtilde(\ustar)$ is the \Frechet\ derivative of $\Jtilde$ at~$\ustar$,
provided that it exists.
In order to justify this procedure, we prove that the control-to-state operator $\calS$ is \Frechet\ differentiable,
so that $D\Jtilde(\ustar)$ \an{can be computed via} the chain rule.
However, \an{evaluating} the \Frechet\ derivative of $\calS$ involves the so-called {\it linearized system},
which depends on an arbitrary increment.
\an{To address this difficulty}, we introduce a proper {\it adjoint system}
that allows us to give \eqref{preNC} a more suitable form.
Namely, we prove that the above inequality \an{can be expressed as}
\Beq
  \intQ \bigl( b_3 \ustar + \oomega \bigr) \cdot (\uu-\ustar)
  + \intQ \Lambdastar \cdot (\uu-\ustar)
  \geq 0
  \qquad \hbox{for every $\uu\in\Uad$}\,,
  \label{betterNC}
\Eeq
where $\oomega$ is the first component of the solution $\soluza$ to the adjoint problem introduced in Section~\ref{ADJOINT}. \juerg{In Section~\ref{SPARSITY}, we exploit the variational inequality
\eqref{betterNC} to derive a sparsity result for the optimal controls if the sparsity functional 
$G$ is given in the special form~\eqref{Gsparse} (see Theorem~\ref{Sparse} below).}

The \an{remainder} of the paper is organized as follows.    
The existence of an optimal control is proved in Section~\ref{CONTROL}.
The following Section~\ref{FRECHET} is devoted to the \an{analysis} of the control-to-state operator,
namely, to its \Frechet\ differentiability.
The necessary condition \eqref{betterNC} is then proved in Section~\ref{ADJOINT},
where a natural adjoint problem is introduced and solved\juerg{, and the sparsity is
addressed in the final Section~\ref{SPARSITY}}.

\juerg{Throughout the paper, we assume that $\eps =m = \eta_0 =1$ without loss of generality  
and, besides the \Holder\ and Cauchy--Schwarz inequalities}
and the Korn inequality \eqref{korn} already mentioned, we \juerg{widely use}  Young's inequality
\Beq
  ab \leq \delta a^2 + \frac 1{4\delta} \, b^2
  \quad \hbox{for every $a,b\in\erre$ and $\delta>0$}.
  \label{young}
\Eeq
We also account for the Sobolev inequalities, 
as well as for some inequalities associated to the elliptic regularity theory 
and to the compact embeddings between Sobolev spaces (via Ehrling's lemma).
In \an{particular}, we have 
\begin{align}
  & \norma z_q
  \leq \CS \, \normaV z
  \quad \hbox{for every $z\in V$ and $q\in[1,6]$},
  \label{sobolev1}
  \\[2mm]
  \separa
  & \norma z_\infty
  \leq \CS \, \normaW z
  \quad \hbox{for every $z\in \an{W}$},
  \label{sobolev2}
  \\[2mm]
  \separa
  & \normaW z
  \leq \CE \, \bigl( \norma{\Delta z} + \norma z \bigr)
  \quad \hbox{for every $z\in W$},
  \label{elliptic1}
  \\[2mm]
  & \pier{\norma z_{\Hx3}
  \leq \CE \, \bigl( \norma{\nabla \Delta z} + \norma z \bigr)
  \quad \hbox{for every $z\in\Hx3\last{{}\cap W{},}$}}
  \label{elliptic1.5}
  \\[2mm]
  & \norma z_{\Hx4}
  \leq \CE \, \bigl( \norma{\Delta^2 z} + \norma z \bigr)
  \quad \hbox{for every $z\in\Hx4$ with $z,\Delta z\in W$},
  \label{elliptic2}
  \\[2mm]
  \separa
  & \normaV z 
  \leq \delta \, \norma{\Delta \pier{{}z{}} }+ \Cdelta \, \norma z  
  \quad \hbox{for every $z\in W$ \pier{and every $\delta>0$}},
  \label{compact1}
  \\[2mm]
    & \pier{ \norma z_{\Hx2}
  \leq \delta \, \norma{\nabla \Delta \pier{{}z{}} }+ \Cdelta \, \norma z  
  \quad \hbox{for every $z\in \Hx3\cap W$ and every $\delta>0$}},
  \label{compact1.5}
  \\[2mm]
  & \norma z_{\Hx3}
  \leq \delta \, \norma{\Delta^2 z} + \Cdelta \, \normaV z  
  \non
  \\
  & \quad \hbox{for every $z\in\Hx4$ with $z,\Delta z\in W$ and every $\delta>0$} \,.
  \label{compact2}
\end{align}
The constants on the \rhs s of \accorpa{sobolev1}{elliptic2} depend only on~$\Omega$, 
while $\Cdelta$ in \accorpa{compact1}{compact2} depends on both $\Omega$ and~$\delta$.
\an{Similar inequalities hold, of course, for vector-and matrix-valued functions.}
In connection with \eqref{elliptic1} and \eqref{elliptic2}, 
we recall that for $z$ to belong to $W$ it suffices that $z\in V$, $\Delta z\in H$, 
and the homogeneous Neumann boundary condition is satisfied in the usual weak sense;
moreover,  $z$ belongs to $\Hx4$ whenever both $z$ and $\Delta z$ belong to~$W$. \pier{Concerning 
the equivalence of norms actually stated by \eqref{elliptic1.5}, we observe that the inequality 
$$\norma z_{\Hx3}  \leq \CE  \bigl( \norma{\Delta z}_{\Hx1} + \norma z \bigr) \quad \hbox{for all } z \in \Hx3\cap W $$ 
follows by elliptic regularity. \last{Furthermore}, since $z\in W$ satisfies  $\iO \Delta z = \iG \dn z =0$, the Poincar\'e--Wirtinger inequality implies that
\(\, 
\| \Delta z \|_{\Hx1} \le C \, \| \nabla \Delta z \|\,
\)
for some constant $C>0$. Combining the above estimates yields \eqref{elliptic1.5}.}

We conclude this section by introducing a convention for the notation of the constants 
that enter the estimates we are going to perform.
The lowercase symbol $c$ denotes a generic constant
that depends only on $\Omega$, $T$, the structure of the system, and an upper bound for the norms of the data.
Notice that the value of $c$ may change from line to line and even within the same line.
\an{Furthermore, whenever a positive constant depends on a specific parameter such as $\delta$, we indicate this dependence by using a subscript and writing $\cdelta$ instead of a general $c$.}
On the contrary, specific constants we want to refer to are denoted by different symbols,
like in \eqref{stability} and~\eqref{sobolev1}, where different characters are~used.


\section{The existence of an optimal control}
\label{CONTROL}
\setcounter{equation}{0}

\an{We now proceed with the proof of Theorem~\ref{Existencecontrol}.}
We denote by $\Lambda$ the infimum of the functional to be minimized.
Then, $\Lambda$~is a nonnegative real number and there exists a sequence $\graffe\un$ in $\Uad$ such that
\Beq
  \Lambda \leq \calJ(\phin,\un) \leq \Lambda + \frac 1n
  \quad \hbox{for all $n\in\enne$,}
  \label{minimizing}
\Eeq
where $\phin$ is the second component of the solution $\soluzn$ to the state system corresponding to~$\un$.
Since $\Uad$ is bounded and closed in $\an{\LL^\infty(Q)}$, we can
\juerg{without loss of generality} assume that
\Beq
  \un \to \ustar
  \quad \hbox{weakly star in $\an{\LL^\infty(Q)}$} 
  \non
\Eeq
for some $\ustar\in\Uad$.
On the other hand, $\soluzn$ satisfies the stability estimate \eqref{stability}
with a constant $K_1$ that does not depend on~$n$.
Therefore, we  \an{may also} assume that\an{, as $n\ \to \infty$,}
\begin{align}
   \vn \to \vstar
  & \quad \hbox{weakly in $\L2\VV$},
  \label{convvn}
  \\
   \phin \to \phistar
  & \quad \hbox{weakly star in $\H1\Vp\cap\L\infty W\cap\L2{\Hx5}$},
  \label{convphin}
  \\
   \mun \to \mustar
  & \quad \hbox{weakly in $\L2V$},
  \label{convmun}
  \\
   \wn \to \wstar
  & \quad \hbox{weakly star in $\L\infty H\cap\L2{\Hx3\cap W}$} \,.
  \label{convwn}
\end{align}
In particular, we have that $\phistar(0)=\phiz$.
Moreover, by applying, e.g., \cite[Sect.~8, Cor.~4]{Simon}, the regularity of~$f$
and well-known compact embeddings,
we deduce the strong convergence properties\an{, as $n\ \to \infty$,}
\begin{align}
  & \phin \to \phistar \,, \quad 
  f(\phin) \to f(\phistar)\,,
  \aand
  f'(\phin) \to f'(\phistar)\,,
  \quad \hbox{uniformly in $\overline Q$},
  \non
  \\
  & \phin \to \phistar
  \quad \hbox{strongly in $\C0{\Wx{1,4}}$}\,,
  \non
\end{align}
whence also
\Beq
  \mun \nabla\phin \to \mustar \nabla\phistar
  \quad \hbox{weakly in $\L2\HH$} \,.
  \non
\Eeq
Therefore, it is \sfw\ to conclude that $\soluzstar$ verifies the time-integrated versions 
of the variational equations \accorpa{prima}{quarta} associated with~$\ustar$\juerg{, for} 
time-dependent test functions.
This means that $\soluzstar$ is the solution to the state system corresponding to~$\ustar$.
\an{Moreover}, the above strong convergence implies that
$\phin(T)$ converges to $\phistar(T)$ strongly in~$H$.
\an{Recalling \eqref{minimizing}, and using the lower semicontinuity of~$\calJ$, we then obtain}
\Beq
  \calJ(\phistar,\ustar)
  \leq \liminf_{n\to\infty} \calJ(\phin,\un)
  \leq \liminf_{n\to\infty} \Bigl( \Lambda + \frac 1n \Bigr)
  = 
  \an{\lim_{n\to\infty} \Bigl( \Lambda + \frac 1n \Bigr)=}
  \Lambda\,,
  \non
\Eeq
which readily yields that $\calJ(\phistar,\ustar)=\Lambda$,
and the proof is complete.


\section{The control-to-state mapping}
\label{FRECHET}
\setcounter{equation}{0}

The {\it control-to-state operator} \an{was} introduced in Section~\ref{STATEMENT} in a \an{preliminary} form.
\an{We now provide its precise definition.} \juerg{To this end, we assume that}
\Beq
\label{defUR}
\hbox{$\UR$ is an open ball in $\L2\HH$ containing $\Uad$}, 
\Eeq 
\juerg{and we introduce the solution space}
\begin{align}
  & \calX
  := \L2\VVz \times
  \bigl(
    \H1\Wp \cap \C0V \cap \L2{\Hx4\cap W}
  \bigr) \times {}
  \non
  \\
  & \qquad {} \times \L2H \times \L2W\,. 
  \label{defX}
  \end{align}
  \juerg{We then define
  \begin{align}
  &\calS:\UR\to\calX\, \pier{{}, \, \hbox{ with }\, } \UR\ni\uu \mapsto (\vv,\phi,\mu,w) =: \calS(\uu) \non \\
  & = \mbox{ the solution to the state system \Pbl\ corresponding to $\uu$} .
  \label{defS}
\end{align}
}
\quad Our proof of the \Frechet\ differentiability of $\calS$ relies on an improvement 
of the continuous dependence inequality~\eqref{contdep}.
We have the following result.
\Blem
\label{Berlin}
Under the same assumptions, and with the  \juerg{notation used in} Theorem~\ref{Contdep}, 
we have that
\Beq
  \norma\phi_{\L\infty W \cap \L2{\Hx5}}
  \leq K_3 \, \norma\uu_{\L2\HH}\,,
  \label{berlin}
\Eeq
\an{where the constant $K_3$ \juerg{has the same dependencies} as the constant $K_2$ in} \eqref{contdep}.
\Elem

\Bdim
We proceed formally, for brevity, since a rigorous proof \an{would be rather \juerg{lengthy}}.
\an{A rigorous approach could follow the procedure used} in~\cite{CGSS8} to prove Theorem~\ref{Wellposedness},
which is based on the approximation of Problem \Pbl\ by a Faedo--Galerkin scheme
constructed by means of the eigenfunctions of related eigenvalue problems.
In particular, as far as the component $\phi$ is concerned, 
the eigenvalue problem for the Laplace operator with homogeneous Neumann boundary conditions has been chosen.
This provides a high regularity to the discrete solutions \juerg{$\phi_n$ in the 
$n$\last{-}dimensional approximating spaces. In particular, $\Delta^m\phi_n$ belongs 
to $W$ for every nonnegative integer $m\in\enne$}.
\an{Since the solution \pier{to Problem \Pbl\ is unique}, the solutions appearing in the statement of the lemma coincide with those obtained through this discretization procedure \pier{and the subsequent passage to the limit as $n\to \infty$}. Consequently, formal estimates performed directly on the solution to Problem \Pbl\ can be justified via the discrete scheme. Accordingly, we confine ourselves to providing a formal proof of \eqref{berlin} by directly testing certain equations with functions related to the solution, whose discrete counterparts would be admissible in the Faedo--Galerkin scheme.}

%

We take the difference of \eqref{seconda} written for both solutions,
and we do the same for \eqref{quinta}.
We have, \aet,~that
\begin{align}
  & \iO \dt\phi \, z 
  + \iO \bigl(
    \vv \cdot \nabla\phi_1 + \vv_2 \cdot \nabla\phi
  \bigr) z
  + \iO \nabla\mu \cdot \nabla z
  = \iO \bigl( S(\phi_1) - S(\phi_2) \bigr) z \,,
  \label{dseconda}
  \\[1mm]
  & 
  \an{\iO \nabla \bigl( 
    - \Delta\phi + f(\phi_1) - f(\phi_2) 
  \bigr) \cdot \nabla  z}
  + \iO \bigl[ 
    \bigl( f'(\phi_1) - f'(\phi_2) \bigr) (-\Delta\phi_1) 
    + f'(\phi_2) (-\Delta\phi) 
  \bigr] z
  \non
  \\
  & \quad {}
  + \iO \bigl(
    f(\phi_1) f'(\phi_1) - f(\phi_2) f'(\phi_2) 
  \bigr) z 
  + \nu \iO \bigl(
    -\Delta\phi + f(\phi_1) - f(\phi_2) 
  \bigr) z
  = \iO \mu z \,,
  \label{dquinta}
\end{align}
\an{where $z\in V$ is arbitrary in \eqref{dseconda} 
and 
\eqref{dquinta}.}
\an{No duality is needed for the time derivative, as it reduces to an integral in the corresponding discrete framework.}
Now, we test \eqref{dseconda} by $\Delta^2\phi$ and \eqref{dquinta} by $-\Delta^3\phi$,
take the sum of the resulting equalities, and rearrange a little.
On the \lhs,  we \an{retain} the terms
\Beq
  \frac 12 \, \ddt \, \norma{\Delta\phi}^2
  + \iO |\nabla\Delta^2\phi|^2 \,.
  \label{lhs}
\Eeq
The sum of the terms involving $\mu$ on the \rhs\ is given~by
\Beq
  - \iO \nabla\mu \cdot \nabla\Delta^2\phi  
  + \iO \mu (-\Delta^3\phi)
  = 0 \,,
  \non
\Eeq
\an{using integration by parts and the homogeneous Neumann boundary conditions.}
Now, we estimate, \aet, the remaining terms on the \rhs.
More precisely, we prepare the estimates that produce terms 
that can be controlled by the \lhs~\eqref{lhs} (either directly or via the Gronwall lemma after time integration) 
and terms whose integrals over $(0,T)$ can be estimated by the \rhs\ of~\eqref{contdep}.
We also make a direct use of the $L^\infty$ in time estimates for the solutions ensured by~\eqref{stability}. 
In particular, we can assume that all of the nonlinearities are \Lip\ continuous and bounded.
We repeatedly account for the Young and \Holder\ inequalities and for some of the inequalities \accorpa{sobolev1}{compact2}.
Finally, $\delta$~is a positive parameter whose value \an{will be} chosen later on.
We \an{then have}
\Beq
  - \iO \bigl(
    \vv \cdot \nabla\phi_1
  \bigr) \Delta^2\phi
  \leq c \, \norma\vv_4 \, \norma{\nabla\phi_1}_4 \, \norma{\Delta^2\phi}
  \leq c \, \normaV\vv^2
  + c \, \norma\phi_{\Hx4}^2\,,
  \non 
\Eeq
where we have also used the continuity of the embedding\an{%
$$\L\infty W \emb \L\infty {\Wx {1,4}}$$
and the fact that} $\phi_1$ is estimated in $\L\infty W$.
Next, we have
\Beq
  - \iO (\vv_2 \cdot \nabla\phi) \Delta^2\phi
  \leq \norma{\vv_2}_4 \, \norma{\nabla\phi} \, \norma{\Delta^2\phi}_4
  \leq \delta \, \normaV{\Delta^2\phi}^2
  + \cdelta \, \normaV{\vv_2}^2 \, \normaV\phi^2 \,.
  \non 
\Eeq
\an{The integral involving $S$ is straightforward to handle and therefore omitted.
\juerg{The next term, however,  requires a} careful treatment. Noting} that in the discrete case $f(\phi_i)$, $\Delta\phi_i$, \pier{$\Delta^2 \phi_i$} and $\Delta f(\phi_i)$ 
satisfy the homogeneous Neumann boundary condition
(for the last one, see below for the computation of its gradient),
the first \pier{two integrations by parts are} justified\an{, yielding}
\begin{align}
  & 
 	\an{- \iO \nabla \bigl( 
    f(\phi_1) - f(\phi_2) 
  \bigr) \cdot \nabla \Delta^3\phi} =
  \pier{\iO   \Delta \bigl(f(\phi_1) - f(\phi_2)  \bigr) \Delta^3\phi}
  \non
   \\
  &{}= \an{-} \iO \nabla\Delta\bigl( 
    f(\phi_1) - f(\phi_2)
  \bigr) \cdot \nabla\Delta^2\phi
  \non
  \\
  & {}\leq \delta \, \norma{\nabla\Delta^2\phi}^2
  + \cdelta \, \norma{\nabla\Delta\bigl( f(\phi_1) - f(\phi_2) \bigr)}^2 \,.
  \non 
\end{align}
The treatment of the last term needs a preliminary computation.
Namely, \pier{we have that}
\begin{align}
  & \Delta\bigl( f(\phi_1) - f(\phi_2) \bigr)
  \non
  \\
  & = f''(\phi_1) |\nabla\phi_1|^2 - f''(\phi_2) |\nabla\phi_2|^2
  + f'(\phi_1) \Delta\phi_1 - f'(\phi_2) \Delta\phi_2\,,
  \non
\end{align}
whence, denoting by $D^2\phi_i$ the Hessian matrix of~$\phi_i$,
\begin{align}
  & \nabla\Delta\bigl( f(\phi_1) - f(\phi_2) \bigr)
  \non
  \\
  & = f\an{^{(3)}}(\phi_1) |\nabla\phi_1|^2 \nabla\phi_1 - f\an{^{(3)}}(\phi_2) |\nabla\phi_2|^2 \, \nabla\phi_2
  + 2f''(\phi_1) D^2\phi_1 \nabla\phi_1 - 2f''(\phi_2) D^2\phi_2 \nabla\phi_2
  \non
  \\
  & \quad {}
  + f''(\phi_1) \Delta\phi_1 \nabla\phi_1 - f''(\phi_2) \Delta\phi_2 \nabla\phi_2
  + f'(\phi_1) \nabla\Delta\phi_1 - f'(\phi_2) \nabla\Delta\phi_2 \,.
  \non
\end{align}
Hence, we \juerg{need to} estimate these four differences.
However, we consider just the first two, since the others can be dealt with in a similar way.

\juerg{At this point, we use the fact that there exists some constant $\hat c>0$ (actually,
a little calculus shows that the minimal such constant is $\hat c= 3/2$) such that}
\begin{align*}
	\an{\bigl||x|^2x-|y|^2y\bigr|\leq \juerg{\hat c}\, |x-y|(|x|^2+|y|^2),
 \quad\mbox{for all }\, x,y\in\erre^3\,.}
\end{align*}
\juerg{We therefore can estimate as follows:} 
\begin{align}
  & \norma{f\an{^{(3)}}(\phi_1) |\nabla\phi_1|^2 \nabla\phi_1 - f\an{^{(3)}}(\phi_2) |\nabla\phi_2|^2 \, \nabla\phi_2}^2
  \non
  \\
  & = \norma{(f\an{^{(3)}}(\phi_1)-f\an{^{(3)}}(\phi_2))|\nabla\phi_1|^2\nabla\phi_1 + f\an{^{(3)}}(\phi_2)(|\nabla\phi_1|^2\nabla\phi_1-|\nabla\phi_2|^2\nabla\phi_2)}^2
  \non
  \\
  & \leq c \, \norma{|\phi|\,|\nabla\phi_1|^3}^2
  + c \, \norma{|\nabla\phi|\,(|\nabla\phi_1|^2+|\nabla\phi_2|^2)}^2
  \non
  \\
  \separa
  & \leq c \, \norma\phi_6^2 \, \norma{\nabla\phi_1}_\infty^2 \, \norma{\nabla\phi_1}_6^4
  + c \, \norma{\nabla\phi}_6^2 \bigl( \norma{\nabla\phi_1}_6^4 + \norma{\nabla\phi_2}_6^4 \bigr) 
  \non
  \\
  \separa
  & \leq c \, \normaV\phi^2 \, \norma{\phi_1}_{\Hx3}^2 \, \normaW{\phi_1}^4
  + c \, \normaW\phi^2 \bigl( \normaW{\phi_1}^4 + \normaW{\phi_2}^4 \bigr) 
  \non
  \\
  & \leq c \, \normaV\phi^2 \, \norma{\phi_1}_{\Hx3}^2 
  + c \, \normaW\phi^2 \,.
  \non
\end{align}
As for the second difference, we have that
\begin{align}
  & \norma{2f''(\phi_1) D^2\phi_1 \nabla\phi_1 - 2f''(\phi_2) D^2\phi_2 \nabla\phi_2}^2
  \non
  \\
  & \leq c \, \norma{\bigl( f''(\phi_1)-f''(\phi_2) \bigr) D^2\phi_1 \nabla\phi_1}^2
  + c \, \norma{f''(\phi_2) \bigl( D^2\phi_1- D^2\phi_2 \bigr) \nabla\phi_1}^2
  \non
  \\
  & \quad {}
  + c \, \norma{f''(\phi_2) D^2\phi_2 \bigl( \nabla\phi_1 - \nabla\phi_2 \bigr)}^2
  \non
  \\
  \separa
  & \leq c \, \norma\phi_6^2 \, \norma{D^2\phi_1}_6^2 \, \norma{\nabla\phi_1}_6^2
  + c \, \norma{D^2\phi}_4^2 \, \norma{\nabla\phi_1}_4^2 
  + c \, \norma{D^2\phi_2}_4^2 \, \norma{\nabla\phi}_4^2
  \non
  \\
  \separa
  & \an{\leq c \, \normaV\phi^2 \, \norma{\phi_1}_{\Hx3}^2\, \norma{\phi_1}_{W}^2
  + c \, \norma\phi_{\Hx3}^2\, \norma{\phi_1}_{W}^2
  + c \, \norma{\phi_2}_{\Hx3}^2 \, \normaW\phi^2 
  \non}
    \\ 
  & \an{\leq c \, \normaV\phi^2 \, \norma{\phi_1}_{\Hx3}^2
  + c \, \norma\phi_{\Hx3}^2
  + c \, \norma{\phi_2}_{\Hx3}^2 \, \normaW\phi^2 \,.
  \non}
\end{align}
Let us come to the next two terms \juerg{which originate from \eqref{dquinta}.} \pier{Concerning the first, we}
have
\begin{align}
  & - \iO \bigl( f'(\phi_1) - f'(\phi_2) \bigr) (-\Delta\phi_1) (-\Delta^3\phi)
  \non
  \\
  & = \iO \nabla \bigl[ (-\Delta\phi_1) \bigl( f'(\phi_1) - f'(\phi_2) \bigr) \bigr] \cdot \nabla\Delta^2\phi
  \non
  \\
  & = \iO (-\Delta\phi_1) \bigl( f''(\phi_1) - f''(\phi_2) \bigr) \nabla\phi_1 \cdot \nabla\Delta^2\phi
  + \iO (-\Delta\phi_1) f''(\phi_2) \nabla\phi \cdot \nabla\Delta^2\phi
  \non
  \\
  & \quad {}
  + \iO \bigl( f'(\phi_1) - f'(\phi_2) \bigr) \nabla(-\Delta\phi_1) \cdot \nabla\Delta^2\phi
  \non
  \\
  \separa
  & \leq \delta \iO |\nabla\Delta^2\phi|^2
  + \cdelta \, \norma{\Delta\phi_1}_6^2 \an{\, \norma\phi_6^2\, \norma{\nabla\phi_1}_6^2 }
  + \cdelta \, \norma{\Delta\phi_1}_4^2 \, \rev{\norma{\nabla\phi}_4^2}
  + \cdelta \, \norma\phi_4^2 \, \norma{\nabla\Delta\phi_1}_4^2 
  \non
  \\
  \separa
  & \leq \delta \iO |\nabla\Delta^2\phi|^2
  \non
  \\
  & \quad {}
  + \cdelta \, \norma{\phi_1}_{\Hx3}^2 \, \normaV\phi^2
  + \cdelta \, \norma{\phi_1}_{\Hx3}^2 \bigl( \norma\phi^2 + \norma{\Delta\phi}^2 \bigr)
  + \cdelta \, \normaV\phi^2 \, \norma{\phi_1}_{\Hx4}^2\,, 
  \non 
\end{align}
\an{where we also have repeatedly used \juerg{that the mapping $t \mapsto \|\phi_1(t)\|_{W}$ belongs to $L^\infty(0,T)$.} Next, we \pier{deduce that}}
\begin{align}
  & - \iO f'(\phi_2) (-\Delta\phi) (-\Delta^3\phi)
  = \iO \nabla \bigl( f'(\phi_2) (-\Delta\phi) \bigr) \cdot \nabla\Delta^2\phi
  \non
  \\
  & \leq \delta \iO |\nabla\Delta^2\phi|^2
  + \cdelta \, \norma{\nabla \bigl( f'(\phi_2) \an{\Delta\phi} \bigr)}^2
  \non
  \\
  & = \delta \iO |\nabla\Delta^2\phi|^2
  + \cdelta \, \norma{\an{\Delta\phi}f''(\phi_2)\nabla\phi_2 + f'(\phi_2)\nabla\Delta\phi}^2
  \non
  \\
  \separa
  & \leq \delta \iO |\nabla\Delta^2\phi|^2
  + \cdelta \,\an{\, \norma{\Delta\phi}^2 \norma{\nabla\phi_2}_\infty^2 }
  + \cdelta \, \norma{\nabla\Delta\phi}^2
  \non
  \\
  & \leq \delta \iO |\nabla\Delta^2\phi|^2
  + \cdelta \an{\, \norma{\Delta\phi}^2\, \norma{\phi_2}_{\Hx3}^2 }
  + \cdelta \, \norma\phi_{\Hx3}^2 \,.
  \non 
\end{align}
As for the next term, we replace $f\!f'$ by a generic smooth function~$g$.
The inequality we obtain, \juerg{when} applied to $g=-\nu f$, \juerg{yields} an estimate 
for the last term \juerg{to be considered}.
We have
\begin{align}
  & - \iO \bigl( g(\phi_1) - g(\phi_2) \bigr) (-\Delta^3\phi) 
  = \iO \nabla \bigl( g(\phi_1) - g(\phi_2) \bigr) \cdot \nabla\Delta^2\phi
  \non
  \\
  & = \iO \bigl( (g'(\phi_1)-g'(\phi_2))\nabla\phi_1 + g'(\phi_2)\nabla\phi \bigr) \cdot \nabla\Delta^2\phi 
  \non
  \\
  \separa
  & \leq \delta \iO |\nabla\Delta^2\phi|^2
  + \cdelta \, \bigl( \norma{\pier{|\phi| \, |\nabla\phi_1|}}^2 + \norma{\nabla\phi}^2 \bigr)
  \non
  \\
  \separa
  & \leq \delta \iO |\nabla\Delta^2\phi|^2
  + \cdelta \, \bigl( \norma\phi_4^2 \, \norma{\nabla\phi_1}_4^2 + \norma{\nabla\phi}^2 \bigr)
  \non
  \\
  & \leq \delta \iO |\nabla\Delta^2\phi|^2
  + \cdelta \, \bigl( \normaW{\phi_1}^2 + 1 \bigr) \normaV\phi^2 
    \non
  \\
  & \leq \an{\delta \iO |\nabla\Delta^2\phi|^2
  + \cdelta  \normaV\phi^2 \,.}
  \non 
\end{align}
Finally, \pier{it is clear that}
\begin{align}
  & \pier{{}- \nu \iO (-\Delta\phi) (-\Delta^3\phi)
  =  \nu} \iO \nabla\Delta\phi \cdot \nabla\Delta^2\phi
  \leq \delta \iO |\nabla\Delta^2\phi|^2
  + \cdelta \, \norma\phi_{\Hx3}^2 \,.
  \non 
\end{align}
At this point, we recall the \lhs\ \eqref{lhs} of the identity obtained  as said at the beginning,
and collect all \juerg{the above} inequalities to estimate the corresponding \rhs.
Then, we choose $\delta$ small enough, rearrange, and integrate over $(0,t)$ for an 
arbitrary $t\in(0,T]$.
\pier{Taking also~\eqref{contdep} into account, we finally obtain that}
\Beq
  \norma{\Delta\phi(t)}^2
  + \intQt |\nabla\Delta^2\phi|^2
  \leq c \, \norma\uu_{\L2\HH}^2
  + \iot \psi(s) \, \norma{\Delta\phi(s)}^2 \, ds\,,
  \non
\Eeq
\an{where $\psi$ is naturally defined by the above estimates. We observe that
\juerg{the most delicate terms arising in $\psi$ involve the expressions} $\|\phi_1\|_{\Hx4}^2$ and $\|\phi_2\|_{\Hx3}^2$. Therefore, in view of \eqref{stability}, \juerg{the
mapping} $t \mapsto \psi(t)$ is bounded in $L^1(0,T)$. 
Integrating in time and applying Gronwall's lemma, we can therefore deduce that}
\Beq
  \norma{\Delta\phi}_{\L\infty H}
  + \norma{\nabla\Delta^2\phi}_{\L2H}
  \leq c \, \norma\uu_{\L2\HH}\,,
  \non
\Eeq
whence \eqref{berlin} follows from \pier{\eqref{contdep}, \eqref{elliptic1}, \eqref{elliptic1.5} and} elliptic regularity.
\Edim
\Brem
\pier{Having proved Lemma~\ref{berlin}, and adopting the same notation 
as in
Theorem~\ref{Contdep}, we observe that, by first taking the difference
of equations~\eqref{P-quarta} and then of~\eqref{P-terza}, and arguing as
in the above proof in order to control the derivatives of the nonlinear
terms, we are led to the further estimate
\begin{align*}
  \norma w_{\L\infty H \cap \L2{\Hx3}} + \norma \mu_{\L2 V}
  \le K_4 \, \norma \uu_{\L2 \HH},
\end{align*}
which holds for the differences $w = w_1 - w_2$, $\mu = \mu_1 - \mu_2$,
and $\uu = \uu_1 - \uu_2$. The constant $K_4$ has the same dependencies
as $K_3$ in~\eqref{berlin}.}
\Erem

From now on, \juerg{we generally assume:} 
\Beq
  \vbox{\hsize .7\hsize \noindent\it
    \juerg{All of the above assumptions on the structure o\rev{f} the original system, 
     the initial datum,  the cost functional, and} the set of  admissible controls, are satisfied.}
  \label{allassumptions}
\Eeq
For the reader's convenience, we recall that these assumptions are given in \Hpstruttura, \eqref{regphiz}, \an{and \eqref{cost}--\eqref{hpuminmax}.}

\an{We now turn \juerg{our interest} to the study of the \Frechet\ \juerg{differentiability}
 of the solution operator~$\calS$. As a preliminary step, we observe that such a derivative is characterized through the associated {\it linearized system}, which we therefore introduce and analyze next.}
For its solution, we can use the same notation \juerg{as that} adopted for the solution to problem \Pbl,
since no confusion can arise. 
Indeed, the solution $\calS(\ustar)$ to the original problem corresponding to a fixed $\ustar\in\UR$
is denoted by \an{$\soluzstar$.}
Given an element $\hh\in\L2\HH$, the associated linearized problem consists in looking for a quadruplet $\soluz$ with the regularity
\begin{align}
  & \vv \in \L2\VVz ,
  \label{regvl}
  \\
  & \phi \in \H1\Wp \cap \C0V \cap \L2{\Hx4\pier{{}\cap W{}}} ,
  \label{regphil}
  \\
  & \mu \in \L2H ,
  \label{regmul}
  \\
  & w \in \L2W,
  \label{regwl}
\end{align}
\Accorpa\Regsoluzl regvl regwl
that solves the variational equations
\begin{align}
  & \iO D\vv : \nabla\zz \pier{{}+ \iO \lambda(\phistar) \vv \cdot \zz}
  \an{{}+ \iO \lambda'(\phistar) \phi \vstar \cdot \zz}
  \non
  \\
    & = \iO \bigl(
    \mu \nabla\phistar
    + \mustar \nabla\phi
    + \hh
  \bigr) \cdot \zz
  \qquad \hbox{for every $\zz\in\VVz$ and \aet},
  \label{primal}
  \\
  \separa
  & \< \dt\phi , z >
  + \iO \vv \cdot \nabla\phistar \, z
  + \iO \vstar \cdot \nabla\phi \, z
  - \iO \mu \, \Delta z
  \non
  \\
  & = \iO S'(\phistar) \phi z
  \qquad \hbox{for every $z\in W$ and \aet},
  \label{secondal}
  \\
  & \iO \nabla w \cdot \nabla z
  + \iO f''(\phistar) \phi \wstar z
  + \iO \an{(f'(\phistar) +\nu)} w z 
  \non
  \\
  & = \iO \mu z
  \qquad \hbox{for every $z\in V$ and \aet},
  \label{terzal}
  \\
  & \an{\iO \nabla \phi\cdot \nabla z}
  + \iO f'(\phistar) \phi z \pier{{} = \iO w z }
  \qquad \hbox{for every $\an{z\in V}$ and \aet}
  \label{quartal}
\end{align}
and satisfies the initial condition
\Beq
  \phi(0) = 0 \,.
  \label{cauchyl}
\Eeq
\Accorpa\Pbll primal cauchyl
Also in this case, one can eliminate the fourth component \an{$w$} of the solution \an{$\soluz$} whenever \juerg{this  is convenient}. \pier{Indeed, \eqref{quartal} and the regularities \eqref{regphil} and \eqref{regwl}
yield
\Beq
- \Delta \phi + f'(\phistar) \phi =  w 
 \last{\quad \aeQ \, .}
  \label{4pierl}
\Eeq
Hence,} if $\soluz$ is a solution to \pier{\accorpa{primal}{quartal},} then the triplet $(\vv,\phi,\mu)$ solves \pier{\eqref{primal}, \eqref{secondal} and} the variational equation
\begin{align}
  & 
\an{ \iO \nabla \bigl(-\Delta\phi + f'(\phistar) \phi
  \bigr) \cdot \nabla z}
  + \iO f''(\phistar) \phi \, \bigl(
    -\Delta\phistar + f(\phistar)
  \bigr) z
  \non
  \\
  & \quad {}
  + \iO \bigl( f'(\phistar) + \nu \bigr) \bigl( -\Delta\phi + f'(\phistar) \phi \bigr) z
  \non
  \\
  & = \iO \mu z
  \qquad \hbox{for every $\an{z\in V}$ and \aet} \,.
  \label{quintal}
\end{align}
Conversely, if the \pier{pair  $(\phi,\mu)$} solves \eqref{quintal} and \eqref{quartal} is taken as a definition of~$w$,
then \eqref{terzal} is satisfied as well.
So, we can consider either of these equivalent problems.
We also notice that the equations \eqref{terzal}, \eqref{quartal} and \eqref{quintal} 
can be written in \an{strong} form \an{as} boundary value problems,
thanks to the regularity of the solution required in \Regsoluzl.
For instance, \an{regarding} \eqref{quintal}, we have
\begin{align}
  & - \Delta \bigl(
    -\Delta\phi + f'(\phistar) \phi
  \bigr) 
  + f''(\phistar) \phi \, \bigl(
    -\Delta\phistar + f(\phistar)
  \bigr) 
  \non
  \\
  & \quad {}
  + \bigl( f'(\phistar) + \nu \bigr) \bigl( -\Delta\phi + f'(\phistar) \phi \bigr) 
  = \mu 
  \qquad \aeQ\,, 
  \label{quintals}
\end{align}
with the boundary conditions
\Beq
  \dn\phi = \dn\Delta\phi = 0 
  \quad \hbox{on $\Sigma$} \,.
  \label{bcquintals}
\Eeq
Recalling the definition of  $\calX$ given in~\eqref{defX}, we have the following result\an{.}

\Bthm
\label{Wellposednessl}
Let $\ustar\in\UR$ be given  and  $\soluzstar:=\calS(\ustar)$ be the corresponding state. 
Then, for every $\hh\in\L2\HH$, the linearized problem \Pbll\ has a unique solution $\soluz$ with the regularity \Regsoluzl.
\rev{Moreover}, the inequality
\Beq
  \norma\soluz_\calX \leq K \, \norma\hh_{\L2\HH}
  \label{stimaperF}
\Eeq
holds true with a constant $K$ that does not depend on~$\hh$, that~is,
\Beq
  \hbox{the linear mapping $\hh \mapsto \soluz$ belongs to $\calL(\L2\HH,\calX)$}.
  \label{perF}
\Eeq

\Ethm

\Bdim
Also in this case, we proceed formally.
However, we make some observations at the end of the proof.
We recall that $\calS(\ustar)$ satisfies the stability estimate~\eqref{stability},
which, in particular, implies that $f(\phistar)$, $f'(\phistar)$ and $f''(\phistar)$ are bounded since $f$ is smooth \an{and $\phistar$ is bounded}. 
Moreover, $\lambda$,~$\lambda'$ and $S'$ are bounded by assumption.
We test \eqref{primal} by~$\vv$ and use the identity~\eqref{identity}.
At the same time, we test \eqref{secondal} by $M(\phi-\Delta\phi)$ 
and \eqref{quintal} 
by $-M\Delta\phi+M\Delta^2\phi-N\mu$,
where $M$ and $N$ are positive parameters whose values will be chosen later.
By rearranging a little, we obtain the identities
\begin{align}
  & \iO |D\vv|^2
  + \iO \lambda(\phistar) |\vv|^2
  \nonumber
  \\
  &= - \iO \lambda'(\phistar) \phi \vstar \cdot \vv
  + \iO \bigl(
    \mu \nabla\phistar
    + \mustar \nabla\phi
    + \hh
  \bigr) \cdot \vv \,,
  \label{testprimal}
  \\[2mm]  
  \separa
  & \frac M2 \, \ddt \, \normaV\phi^2
  \an{{}- M \iO \mu \Delta(\phi-\Delta\phi)}
  \non
  \\
  & = - M \iO \vv \cdot \nabla\phistar \, (\phi-\Delta\phi)
  - M \iO \vstar \cdot \nabla\phi \, (\phi-\Delta\phi)
  \non
  \\
  & \quad {} + \an{M } \iO S'(\phistar) \phi (\phi-\Delta\phi) \,,
  \label{testsecondal}
  \\[2mm]
  \separa  
  & M \iO |\nabla\Delta\phi|^2
  + M \iO |\Delta^2\phi|^2
  + N \iO |\mu|^2 \pier{{}- M \iO \mu ( -\Delta\phi + \Delta^2 \phi )}
  \non
  \\
  & = 
  \pier{{} - N \iO \Delta^2\phi \, \mu
  - \iO (-\Delta) \bigl(
    f'(\phistar) \phi
  \bigr) \bigl( -M\Delta\phi + M\Delta^2\phi - N\mu \bigr)}
  \non
  \\
  & \quad {}
  - \iO f''(\phistar) \phi \, \bigl(
    -\Delta\phistar + f(\phistar)
  \bigr) \bigl( -M\Delta\phi + M\Delta^2\phi - N\mu \bigr)
  \non
  \\
  & \quad {}
  - \iO \bigl( f'(\phistar) + \nu \bigr) \bigl( -\Delta\phi + f'(\phistar) \phi \bigr) \bigl( -M\Delta\phi + M\Delta^2\phi - N\mu \bigr) \,.
  \label{testquintal}
\end{align}
Now, we sum up and notice that two terms involving $\mu$ cancel each other.
Then, by applying the coerciveness inequality \eqref{coercive},
we obtain an inequality whose \lhs\ is given~by \juerg{the expression}
\Beq
  \alpha \, \normaV\vv^2
  + \frac M2 \, \ddt \, \normaV\phi^2
  + M \iO |\nabla\Delta\phi|^2
  + M \iO |\Delta^2\phi|^2
  + N \iO |\mu|^2\,.
  \label{lhsl}
\Eeq
We now estimate the terms of the resulting \rhs\ individually.
\an{More precisely, we derive the estimates so that every term can be controlled by the left-hand side of~\eqref{lhsl}, either directly or, after integrating in time, via  Gronwall's lemma, since the coefficients are bounded in $L^1(0,T)$.}
We make a wide use of the Young and \Holder\ inequalities,
as well as of some of the inequalities \accorpa{sobolev1}{compact2},
even with the precise values of the constants that appear there.
At first, we have
\begin{align*}
  & \an{{}  - \iO \lambda'(\phistar) \phi \vstar \cdot \vv
  \pier{{}\leq c \norma \phi \norma{\vstar}_4 \norma{\vv}_4}
  \leq \frac \alpha 8 \, \normaV\vv^2
  + c \, \normaV\vstar^2 \, \norma\phi^2 \,.
  }
\end{align*}
Next, \pier{we deduce that}
\begin{align}
  & \iO \mu \nabla\phistar \cdot \vv
  \leq \norma\mu \, \norma{\nabla\phistar}_4 \, \norma\vv_4
  \leq \an{ \CS }\norma\mu \, \norma{\nabla\phistar}_4 \,\normaV\vv
  \non
  \\
  & \leq \frac \alpha 8 \, \normaV\vv^2
  + \frac {2\CS^2} \alpha \, \norma{\nabla\phistar}_4^2 \, \norma\mu^2
 \, \an{\leq} \, \frac \alpha 8 \, \normaV\vv^2
  + \hat C \, \norma\mu^2\,,
  \non 
\end{align}
where we have set
\Beq
  \hat C := \frac {2\CS^2} \alpha \, \norma{\nabla\phistar}_{\L\infty{\Lx4}}^2 \,,
  \label{hatC}
\Eeq
\an{recalling} that $\phistar\in\L\infty W$ and that $W\emb\Wx{1,4}$.
The next two terms \pier{can be estimated as follows:}
\begin{align}
  & \iO \mustar \nabla\phi \cdot \vv
  + \iO \hh \cdot \vv
  \leq \norma\mustar_4 \, \norma{\nabla\phi} \, \norma\vv_4
  + \norma\hh \, \norma\vv
  \non
  \\
  & \leq \frac \alpha 8 \, \normaV\vv^2
  + c \, \normaV\mustar^2 \, \normaV\phi^2
  + c \, \norma\hh^2 \,.
  \non 
\end{align}
Let us come to the terms that appear on the \rhs\ of \eqref{testsecondal}.
We have
\begin{align}
  & - M \iO \vv \cdot \nabla\phistar \, (\phi-\Delta\phi)
  \leq M \, \norma\vv_4 \, \norma{\nabla\phistar}_4 \, \norma{\phi-\Delta\phi}
    \non
  \\
  &
  \leq \frac \alpha 8 \, \pier{\normaV\vv^2 
+ c \, M^2} 
\, \norma{\nabla\phistar}_4^2 \, \normaW\phi^2
 \leq \frac \alpha 8 \, \normaV\vv^2
  + \pier{  c \, M^2} 
  \, \normaW\phi^2
   \non
  \\
  &\leq \frac \alpha 8 \, \normaV\vv^2
  + \frac M {\juerg{16}} \, \norma{\Delta^2\phi}^2
  + \cM \, \normaV\phi^2 \,,
  \non 
\end{align}
and, similarly,
\begin{align}
  & - M \iO \vstar \cdot \nabla\phi \, (\phi-\Delta\phi)
  \leq M \, \norma\vstar_4 \, \norma{\nabla\phi} \, \norma{\phi-\Delta\phi}_4
  \non
  \\
  & \leq \an{c}\,
\pier{ M }
\, \norma\phi_{\Hx3}^2
  + \pier{  c \, M }
\, \normaV\vstar^2 \, \normaV\phi^2 
  \leq \frac M {\juerg{16}} \, \norma{\Delta^2\phi}^2 
  + \cM \, \normaV\phi^2
  + \pier{  c \, M }
\, \normaV\vstar^2 \, \normaV\phi^2 \,.
  \non 
\end{align}
Finally, we easily see that
\Beq
  M \iO S'(\phistar) \phi (\phi-\Delta\phi)
  \leq  \pier{ c \, M} 
\, \norma\phi \, \norma{\phi-\Delta\phi}
  \leq \frac M {\juerg{16}}\, \norma{\Delta^2\phi}^2
  + \cM \, \normaV\phi^2 \,.
  \non 
\Eeq
\an{The \pier{first term on the} right-hand side term of \eqref{testquintal} is handled similarly, using Young's inequality to obtain}
\Beq
  - N \iO \Delta^2\phi \, \mu
  \leq \frac N8 \, \norma\mu^2
  + 2N \, \norma{\Delta^2\phi}^2 \,.
  \non 
\Eeq

\an{The last three terms, however, require a more careful analysis.}
In order to minimize the calculation, we introduce a generic function~$z\in\L2H$
and a positive parameter~$\kappa$ \juerg{that are specified later}.
We have
\begin{align}
  & - \iO (-\Delta) \bigl( f'(\phistar) \phi \bigr) z
  \pier{{} = \iO \div \bigl( f''(\phistar )\last{\nabla \phistar }\phi + f'(\phistar) \nabla \phi\bigr)z}
  \non
  \\
  & = \iO \bigl(
    f\an{^{(3)}}(\phistar) |\nabla\phistar|^2 \phi
    + f\an{''}(\phistar) \Delta\phistar \, \phi
    + 2 f''(\phistar) \nabla\phistar \cdot \nabla\phi
    + f'(\phistar) \Delta\phi
  \bigr) z
  \non
  \\
  & \leq  \an{c \bigl( \norma{\nabla\phistar}_6^2 \, \norma\phi_6 
    +  \norma{\Delta\phistar}_4 \, \norma\phi_4 
    + \norma{\nabla\phistar}_4 \, \norma{\nabla\phi}_4 
    + \norma{\Delta\phi} \bigr) \norma z}
  \non
  \\
  & \leq \frac \kappa {32} \, \norma z^2
  + \ck \bigl(
    \normaW\phistar^\rev{4} \, \normaV\phi^2
    + \norma\phistar_{\Hx3}^2 \, \normaV\phi^2
    + \normaW\phi^2
  \bigr)
  \non
  \\
  & \leq \frac \kappa {32} \, \norma z^2
  + \ck \bigl(
    \normaW\phi^2
    + \norma\phistar_{\Hx3}^2 \, \normaV\phi^2
  \bigr)
  \non
  \\
  & \leq \frac \kappa {32} \, \norma z^2
  + \frac M8 \, \norma{\Delta^2\phi}^2
  + \ckM \, \normaV\phi^2 
  + \ck \, \norma\phistar_{\Hx3}^2 \, \normaV\phi^2 
  \non
  \\
  & \leq \frac \kappa {32} \, \norma z^2
  + \frac M8 \, \norma{\Delta^2\phi}^2
  + \ckM \bigl( \norma\phistar_{\Hx3}^2 + 1 \bigr) \normaV\phi^2 \,.
  \non
\end{align}
Next, we find that
\begin{align}
  & - \iO f''(\phistar) \phi \, \bigl(
    -\Delta\phistar + f(\phistar)
  \bigr) z
  \leq c \, \norma\phi_4 \, \bigl( \norma{\Delta\phistar}_4 + 1 \bigr) \, \norma z
  \non
  \\
  & \leq c \, \normaV\phi \, \bigl( \norma\phistar_{\Hx3} + 1 \bigr) \, \norma z
  \leq \frac \kappa {32} \, \norma z^2
  + \ck \, \bigl( \norma\phistar_{\Hx3}^2 + 1 \bigr) \,\normaV\phi^2 \,.
  \non
\end{align}
Finally, \pier{it turns out that}
\begin{align}
  & - \iO \bigl( f'(\phistar) + \nu \bigr) \bigl( -\Delta\phi + f'(\phistar) \phi \bigr)  z
  \leq c \, \normaW\phi \, \norma z
  \non
  \\
  & \leq \frac \kappa {32} \, \norma z^2
  + \ck \, \normaW\phi^2
  \leq \frac \kappa {32} \, \norma z^2
  + \frac M8 \, \norma{\Delta^2\phi}^2
  + \ckM \, \normaV\phi^2 \,.
  \non
\end{align}
At this point, we make a suitable choice of $z$ and $\kappa$ in each of the three last estimates.
Namely, we choose first $z=M(-\Delta\phi+\Delta^2\phi)$ and $\juerg{\kappa=1/M}$, 
and then $z=-N\mu$ and $\juerg{\kappa=1/N}$,
and add the two resulting inequalities.
Moreover, we observe that
\Beq
  \norma{-\Delta\phi+\Delta^2\phi}^2
  \leq 2 \, \norma{\Delta\phi}^2 + 2 \, \norma{\Delta^2\phi}^2
  \leq 3 \, \norma{\Delta^2\phi}^2 + c \, \normaV\phi^2 \,.
  \non
\Eeq
We \juerg{then} obtain the following three corresponding estimates:
\begin{align}
  & - \iO (-\Delta) \bigl(
    f'(\phistar) \phi
  \bigr) \bigl( -M\Delta\phi + M\Delta^2\phi - N\mu \bigr)
  \non
  \\
  &\quad \leq \frac M {32} \, \norma{-\Delta\phi+\Delta^2\phi}^2
  + \frac N {32} \, \norma\mu^2
  + \frac M{\an{4}} \, \norma{\Delta^2\phi}^2 
  + \cMN \bigl( \norma\phistar_{\Hx3}^2 + 1 \bigr) \normaV\phi^2
  \non
  \\
  & \quad\leq \frac {\an{11}M} {32} \, \norma{\Delta^2\phi}^2 
  + \frac N {32} \, \norma\mu^2
  + \cMN \bigl( \norma\phistar_{\Hx3}^2 + 1 \bigr) \normaV\phi^2 \,,
  \non 
  \\[2mm]
  & - \iO f''(\phistar) \phi \, \bigl(
    -\Delta\phistar + f(\phistar)
  \bigr) \bigl( -M\Delta\phi + M\Delta^2\phi - N\mu \bigr)
  \non
  \\
  &\quad \leq \frac {3M} {32} \, \norma{\Delta^2\phi}^2
  + \frac N {32} \, \norma\mu^2
  + \cMN \, \bigl( \norma\phistar_{\Hx3}^2 + 1 \bigr) \,\normaV\phi^2 \,,
  \non 
  \\[2mm]
  & - \iO \bigl( f'(\phistar) + \nu \bigr) \bigl( -\Delta\phi + f'(\phistar) \phi
  \bigr) \bigl( -M\Delta\phi + M\Delta^2\phi - N\mu \bigr)
  \non
  \\
  & \quad\leq \frac {7M} {32} \, \norma{\Delta^2\phi}^2 
  + \frac N {32} \, \norma\mu^2
  + \an{\cMN} \, \normaV\phi^2 \,.
  \non 
\end{align}
These are the last estimates that were needed regarding the \rhs\ of the inequality
obtained by adding \accorpa{testprimal}{testquintal} to each other. 
Then, we deduce an estimate for the whole \rhs\ by just adding the above inequalities to each other. 
Since the related \lhs\ is given by~\eqref{lhsl},
we \juerg{arrive at the estimate}
\begin{align}
  & \frac \alpha 2 \, \normaV\vv^2
  + \frac M2 \, \ddt \, \normaV\phi^2
  + M \iO |\nabla\Delta\phi|^2
  \non
  \\
  & \quad {}
  + \Bigl( \frac {\pier{8}M}{32} - 2N \Bigr) \iO |\Delta^2\phi|^2
  + \Bigl( \frac {\an{25}N}{32} - \hat C \Bigr) \iO |\mu|^2
  \non
  \\
  & \leq \cMN \bigl(
    \normaV\vstar^2 + \normaV\mustar^2 + \norma\phistar_{\Hx3}^2 + 1
  \bigr) \normaV\phi^2
  + c \, \norma\hh^2 \,.
  \non
\end{align}
\an{Now we recall that $\hat C$ is fixed and defined in \eqref{hatC}. We therefore can first choose $N$, and subsequently $M$, in such way that all of the coefficients on the left-hand side become positive.}
\an{We then integrate in time and apply the Gronwall lemma, observing again that the terms containing norms of $\soluzstar$ on the right-hand side are time-integrable due to \eqref{stability}.}
This yields that
\begin{align}
  & \norma\vv_{\L2\HH}
  + \norma\phi_{\L\infty V}
  + \norma{\nabla\Delta\phi}_{\L2H}
  + \norma{\Delta^2\phi}_{\L2H}
  \non
  \\
  & \quad {}
  + \norma\mu_{\L2H}
  \leq c \, \norma\hh_{\L2\HH}\,, 
  \non
\end{align}
and we immediately deduce that
\begin{align}
  & \norma\vv_{\L2\HH}
  + \norma\phi_{\L\infty V\cap\L2{\Hx4\pier{{}\cap W{}}}}
  + \norma\mu_{\L2H}
  \leq c \, \norma\hh_{\L2\HH} \,.
  \non
\end{align}
This estimate \juerg{is rigorous when} performed on an arbitrary solution $\soluz$,
 and it implies uniqueness.
Indeed, by linearity, one can assume $\hh=\0$ and deduce that  the three first components of the solution vanish.
Then, $w$~vanishes as well, due to \eqref{quartal}.

\an{Concerning the existence of a solution and the estimate \eqref{stimaperF}, the formal procedure leading to the above bound can be rigorously justified by performing the same computations on a suitable Faedo–Galerkin approximation (in close analogy with the proof of \cite[Thm.~2.1]{CGSS8}). \juerg{Passage} to the limit then yields the result for the actual solution.}
The same is true for the estimate
\Beq
  \norma{\dt\phi}_{\L2\Wp}
  \leq c \, \bigl(
    \norma\vv_{\L2\VV} + \norma\phi_{\pier{\L\infty V}} + \norma\mu_{\L2H}
  \bigr)\,,
  \non
\Eeq
which can be derived by (the discrete version of) \eqref{secondal}.
Then, \eqref{stimaperF} with $\calX$ given by \eqref{defX} 
follows by combining the above inequalities and taking \eqref{quartal} as a 
definition of~$w$.
\Edim

\juerg{At this point, we recall that $\UR$ is endowed with the topology of $\L2\HH$.}
Here is the main result of this section.

\Bthm
\label{Frechet}
The control-to-state operator $\calS:\UR\to\calX$ is \Frechet\ differentiable in $\UR$. \juerg{More precisely, for every} $\ustar\in\UR$, the Frechet\ derivative $D\calS(\ustar)\in\calL(\L2\HH,\calX)$ acts \juerg{as} follows:
if $\hh\in\L2\HH$, then the value $D\calS(\ustar)[\hh]$ is the solution $\soluz$ to the linearized problem \Pbll,
where $\soluzstar:=\calS(\ustar)$.
\Ethm
\Bdim
In order to keep the paper at a reasonable length,
we do not provide a complete proof.
\an{ Instead, we briefly outline the argument and highlight the most delicate steps.}
Let $\ustar\in\UR$ be fixed. \juerg{In the following, we consider increments} 
$\hh\in\L2\HH$ \juerg{having a sufficiently small norm such that 
 $\ustar+\hh\in \UR$, and we denote by $c$ positive constants that do not depend 
 on the special choice of such increments.}
Then the corresponding states $\soluzstar:=\calS(\ustar)$ and $\soluzh:=\calS(\ustar+\hh)$\an{,}
as well as any convex combination \an{thereof}\an{,} satisfy the stability estimate~\eqref{stability}.
Finally, we consider the solution $\soluz$ to the linearized system \Pbll\ corresponding to $\ustar$ and~$\hh$
and set, for convenience,
\Beq
  \soluzF := \soluzh - \soluzstar - \soluz \,.
  \label{notazF}
\Eeq
\juerg{At this point, we note that $\soluz$,  and therefore also $\soluzF$, depends on $\hh$.
In order not to overload the exposition with indi\last{c}es, we have chosen not to indicate this 
dependence in our notation throughout this proof. Also for other quantities like 
$\Lambda_i$ or $R_i$, which will be introduced below, we will suppress the dependence on
$\hh$ in the notation.} 

\juerg{Returning to the proof, we recall \eqref{perF}, which implies that the 
assertion will be  completely proved once we can show that}
\Beq
  \norma{\soluzF}_\calX 
  \leq c \, \norma\hh_{\L2H}^2
  \label{frechet}
\Eeq
with a constant $c$ that does not depend on~$\hh$.
So, we sketch the proof of \eqref{frechet}.
To this end, we \juerg{consider} the problem solved by $\soluzF$,
which is obtained by subtracting both the systems solved by $\soluzstar$ and $\soluz$
\juerg{from} the one solved by~$\soluzh$.
We have
\Beq
  \soluzF \in \calX\,,
  \label{regsoluzF}
\Eeq
and, with the notation given below, the variational equalities
\begin{align}
  & \iO D\xxi : \nabla\zz
  + \iO \lambda(\phistar) \xxi \cdot \zz
  = \iO \LLamu \cdot \zz
  + \iO \LLamd \cdot \zz
  \non
  \\
  & \quad \hbox{for every $\zz\in\VVz$ and \aet},
  \label{primaF}
  \\
  \separa
  & \< \dt\psi , z >
  - \iO \eta \, \Delta z
  = \iO \Lambda_3 z
  + \iO \Lambda_4 z
  \non
  \\
  & \quad \hbox{for every $z\in W$ and \aet},
  \label{secondaF}
  \\
  \separa
  & \iO \nabla\omega \cdot \nabla z
  + \nu \iO \omega z 
  = \iO \eta z 
  + \iO \Lambda_5 z
  \non
  \\
  & \quad \hbox{for every $z\in V$ and \aet},
  \label{terzaF}
  \\
  \separa
  & \iO \nabla\psi \cdot \nabla z
  = \iO \omega z 
  + \iO \Lambda_6
 z
  \non
  \\
  & \quad \hbox{for every $z\in V$ and \aet},
  \label{quartaF}
\end{align}
are satisfied, as well as the initial condition
\Beq
  \psi(0) = 0 \,.
  \label{cauchyF}
\Eeq
\Accorpa\PblF primaF cauchyF
In the above equations, we have introduced the abbreviating notation 
\begin{align}
   \LLamu
  :=&\,\, 
	\an{- \bigl( \lambda(\phih) - \lambda(\phistar) \bigr) (\vh-\vstar)}
  \non
  \\
  & \,
  - \bigl( \lambda(\phih) - \lambda(\phistar) - \lambda'(\phistar) \phi \bigr) \vstar\,,
  \label{defL1}
  \\
  \LLamd
  := &\,\,\mustar \nabla\psi
  + (\muh-\mustar) \nabla(\phih-\phistar)
  + \eta \nabla\phistar ,
  \label{defL2}
  \\
  \Lambda_3
  :=&\, \pier{{}-{}}\vstar \cdot \nabla\psi
  \pier{{}-{}} (\vh-\vstar) \cdot \nabla(\phih-\phistar)
  \pier{{}-{}} \xxi \cdot \nabla\phistar ,
  \label{defL3}
  \\
  \Lambda_4
  :=& \,\,S(\phih)
  - S(\phistar)
  - S'(\phistar) \phi \,,
  \label{defL4}
  \\
  \Lambda_5 
  :=& \,- \bigl( f'(\phih) - f'(\phistar) - f''(\phistar) \phi \bigr) \wstar
  \non
  \\
  & \,
  - \bigl( f'(\phih) - f'(\phistar) \bigr) (\wh-\wstar)
  - f'(\phistar) \an{\omega} \,,
  \label{defL5}
  \\
  \Lambda_6
  :=&\, - \bigl( f(\phih) - f(\phistar) - f'(\phistar) \phi \bigr) \,.
  \label{defL6}
\end{align}
Clearly, due to the regularity of the solution, 
some of the above equations can be written in their strong form.
This is the case for~\eqref{quartaF}, so that $\omega$ can be explicitly written in terms of $\psi$ and~$\Lambda_6$.
Since it turns out that $\Delta\Lambda_6$ belongs to $\L2H$ (see the computation below),
we obtain after a substitution and an integration by parts in \eqref{terzaF} that
\begin{align}
  & \iO (-\Delta)(-\Delta\psi - \Lambda_6) z
  + \nu \iO ( -\Delta\psi - \Lambda_6) z
  = \iO \eta z
  + \iO \Lambda_5 \, z
  \non
  \\
  & \quad \hbox{for every $z\in H$ and \aet} \,.
  \label{quintaF}
\end{align} 
In performing the estimates that are needed to prove \eqref{frechet},
it is convenient to 
\an{rewrite certain differences in suitable forms.
By applying
standard Taylor expansions \juerg{with integral remainder}, 
we easily find  that}
\begin{align}
   \lambda(\phih) - \lambda(\phistar)
  &= \lambda'(\phistar) (\phih-\phistar)
  + R_1 \,,
  \non
  \\
   S(\phih) - S(\phistar)
  &= S'(\phistar) (\phih-\phistar)
  + R_2 \,,
  \non
  \\
   f'(\phih) - f'(\phistar)
  &= f''(\phistar) (\phih-\phistar)
  + R_3 \,,
  \non
  \\
   f(\phih) - f(\phistar)
  &= f'(\phistar) (\phih-\phistar)
  + R_4\,,
  \non
\end{align}
\juerg{where the remainders satisfy, a.e. in $Q$,}
\Beq
  |R_1| + |R_2| + |R_3| + |R_4|
  \leq c \, |\phih-\phistar|^2 \,.
  \label{remainders}
\Eeq
On the other hand, it results that
\begin{align}
  & \LLamu
  = - \bigl( \lambda(\phih) - \lambda(\phistar) \bigr) (\vh-\vstar)
  - \lambda'(\phistar) \psi \vstar
  - R_1 \, \vstar ,
  \label{newL1}
  \\
  & \Lambda_4
  = S'(\phistar) \psi
  + R_2 \,,
  \label{newL4}
  \\
  & \Lambda_5
  = - \bigl( f''(\phistar) \psi + R_3 \bigr) \wstar - f'(\phistar) ( - \Delta\psi - \Lambda_6)
  \non
  \\
  & \qquad \,{}
  - \bigl( f'(\phih) - f'(\phistar) \bigr) 
  \bigl( - \Delta(\phih-\phistar) + f(\phih) - f(\phistar) \bigr) ,
  \label{newL5}
  \\
  & \Lambda_6 
  = - f'(\phistar) \psi - R_4 \,.
  \label{newL6}
\end{align}

At this point, we start estimating.
As in the proof of Theorem~\ref{Wellposednessl}, it is convenient to introduce two positive parameters $M$ and~$N$.
We test \eqref{primaF} by~$\xxi$ and \eqref{secondaF} by $M(\psi-\Delta\psi)$.
At the same time, we test \eqref{quintaF} by $-N\eta$, $-M\Delta\psi$, and $M\Delta^2\psi$. 
By also recalling the identity \eqref{identity}, we obtain
\begin{align}
  & \iO |D\xxi|^2
  + \iO \lambda(\phistar) |\xxi|^2
  = \iO \LLamu \cdot \xxi
  + \iO \LLamd \cdot \xxi \,,
  \label{testprimaF}
  \\[1mm]
  & 
	\an{\frac M2 \,\frac{d}{dt}\, \normaV \psi^2}  
  - M \iO \eta \, \Delta (\psi-\Delta\psi)
  \non
  \\
  &  \quad= M \iO \Lambda_3 (\psi-\Delta\psi)
  + M \iO \Lambda_4 (\psi-\Delta\psi) ,
  \label{testsecondaF}
  \\
  & - N \iO (-\Delta)(-\Delta\psi - \Lambda_6) \eta
  - N \nu \iO ( -\Delta\psi - \Lambda_6) \eta
  \non
  \\
  &\quad = - N \iO |\eta|^2
  - N \iO \Lambda_5 \, \eta \,,
  \label{testquintaF1}
  \\
  & - M \iO (-\Delta)(-\Delta\psi - \Lambda_6) \Delta\psi
  - M \nu \iO ( -\Delta\psi - \Lambda_6) \Delta\psi
  \non
  \\
  & \quad = - M \iO \eta \Delta\psi
  - M \iO \Lambda_5 \, \Delta\psi \,,
  \label{testquintaF2}
  \\
  & M \iO (-\Delta)(-\Delta\psi - \Lambda_6) \Delta^2\psi
  + M \nu \iO ( -\Delta\psi - \Lambda_6) \Delta^2\psi
  \non
  \\
  & \quad = M \iO \eta \Delta^2\psi
  + M \iO \Lambda_5 \, \Delta^2\psi \,.
  \label{testquintaF3}
\end{align}
\juerg{Next, we add  these four identities, noticing that four of the terms
 involving both $\eta$ and 
$\psi$ cancel out. Then, rearranging terms  and invoking the coerciveness inequality \eqref{coercive}, we   obtain an inequality whose \lhs\ is given by the expression}
\Beq
  \alpha \, \normaV\xxi^2
	\an{\,+\,\frac M2 \,\frac{d}{dt} \,\normaV \psi^2}   
  + N \iO |\eta|^2
  + M \iO |\nabla\Delta\psi|^2
  + M \iO |\Delta^2\psi|^2 \,.
  \label{lhsF}
\Eeq

\juerg{We now have to } estimate the terms on the resulting \rhs\
by using the definitions \accorpa{defL2}{defL3} of $\LLamd$ and $\Lambda_3$ 
and the new expressions \accorpa{newL1}{newL6} for the other $\Lambda_i$,
as well as the estimate \eqref{remainders} of the remainders.
However, the techniques that are needed are quite similar to those employed in the proof of Theorem~\ref{Wellposednessl}, i.e.,
a~wide use of the \Holder, Sobolev and Young inequalities and of the compactness inequality \eqref{compact2}.
For this reason, we provide the details only for the most delicate terms, namely,
those involving the Laplacian of~$\Lambda_6$ (see \accorpa{testquintaF1}{testquintaF3}),
which, by \eqref{newL6},  is the sum of two contributions.
We have
\Beq
  \Delta f'(\phistar)
  = f\an{^{(3)}}(\phistar) |\nabla\phistar|^2
  + f'(\phistar) \Delta\phistar\,,
  \non
\Eeq
whence
\Beq
  \Delta(f'(\phistar) \psi)
  = f\an{^{(3)}}(\phistar) |\nabla\phistar|^2 \psi
  + f'(\phistar) \Delta\phistar \, \psi
  + 2 f''(\phistar) \nabla\phistar \cdot \nabla\psi
  + f'(\phistar) \Delta\psi \,.
  \non
\Eeq
\an{To handle the remaining contribution, we use the explicit representation of~$R_4$, which is given by}
\Beq
  R_4 = \Phi \, |\phih-\phistar|^2 \,,
  \quad \hbox{with the notation} \quad
  \Phi := \int_0^1 (1-s) f''(\phistar+s(\phih-\phistar)) \, ds \,.
  \non
\Eeq
It results that
\begin{align}
  & \nabla\Phi
  = \int_0^1 (1-s) f\an{^{(3)}}(\phistar+s(\phih-\phistar)) \, [\nabla\phistar + s \nabla(\phih-\phistar)] \, ds\,,
  \quad \hbox{whence}
  \non
  \\
  & \Delta\Phi
  = \int_0^1 (1-s) f^{(4)}(\phistar+s(\phih-\phistar)) \, |\nabla\phistar + s \nabla(\phih-\phistar)|^2 \, ds
  \non
  \\
  & \qquad\quad {}
  + \int_0^1 (1-s) f\an{^{(3)}}(\phistar+s(\phih-\phistar)) [\Delta\phistar + s \Delta(\phih-\phistar)] \, ds \,.
  \non
\end{align}
Therefore, we derive that
\begin{align}
  & \Delta R_4
  = \Delta\Phi \, |\phih-\phistar|^2
  + 2 \, \nabla\Phi \cdot \nabla (|\phih-\phistar|^2)
  + \Phi \, \Delta (|\phih-\phistar|^2)
  \non
  \\  
  & = \int_0^1 (1-s) f^{(4)}(\phistar+s(\phih-\phistar)) \, |\nabla\phistar + s \nabla(\phih-\phistar)|^2 \, ds \, |\phih-\phistar|^2
  \non
  \\
  & \quad {}
  + \int_0^1 (1-s) f\an{^{(3)}}(\phistar+s(\phih-\phistar)) \, [\Delta\phistar + s \Delta(\phih-\phistar)] \, ds \, |\phih-\phistar|^2
  \non
  \\
  & \quad {}
  + 2 \int_0^1 (1-s) f\an{^{(3)}}(\phistar+s(\phih-\phistar)) [\nabla\phistar + s \nabla(\phih-\phistar)] \, ds \cdot [2 (\phih-\phistar) \nabla(\phih-\phistar)]
  \non
  \\
  & \quad {}
  + \int_0^1 (1-s) f''(\phistar+s(\phih-\phistar)) \, ds \, [\, 2\, |\nabla(\phih-\phistar)|^2 + 2\, (\phih-\phistar) \Delta(\phih-\phistar)\,] \,.
  \non
\end{align}
\an{We are now in a position to estimate the terms in \accorpa{testquintaF1}{testquintaF3} involving $\Delta \Lambda_6$.  
Using the previous rearrangement, their sum can be written as}
\begin{align}
  & \iO (\Delta\Lambda_6) \bigl( N\eta + M \Delta\psi - M \Delta^2\psi \bigr)
  \non
  \\
  & \leq \Bigl( \frac N2 \, \norma\eta^2 + \frac M2 \, \norma{\Delta^2\psi}^2 + \cM \, \normaV\psi^2 \Bigr)
  + \cMN \, \norma{\Delta\Lambda_6}^2\,,
  \non
\end{align}
and it remains to estimate the last norm.
We have
\Beq
  \norma{\Delta\Lambda_6}
  \leq \norma{\Delta(f'(\phistar)\psi)}
  + \norma{\Delta R_4} \,.
  \non
\Eeq
\juerg{At first, we see} that
\Beq
  \norma{\Delta(f'(\phistar)\psi)}
  \leq c \, \norma{\nabla\phistar}_6^2 \, \norma\psi_6
  + c \, \norma{\Delta\phistar}_4 \, \norma\psi_4
  + c \, \norma{\nabla\phistar}_4 \, \an{\norma{\nabla \psi}_4}
  + c \, \norma{\Delta\psi}\,,
  \non
\Eeq
and the last \an{two terms} can be dealt with by invoking the compactness inequality \eqref{compact2}.
On the other hand, we have \aeQ\ that
\begin{align}
  |\Delta R_4|
  &\leq c \, (|\nabla\phih|^2 + |\nabla\phistar|^2) |\phih-\phistar|^2
  + c \, (|\Delta\phih| + |\Delta\phistar|) |\phih-\phistar|^2
  \non
  \\
  & \quad {}
  + c \, (|\nabla\phih| + |\nabla\phistar|) |\phih-\phistar| \, |\nabla(\phih-\phistar)|
  \non
  \\
  & \quad {}
  + c \, |\nabla(\phih-\phistar)|^2
  + c \, |\phih-\phistar| \, |\Delta(\phih-\phistar)|\,,
  \non
\end{align}
whence
\begin{align}
   \norma{\Delta R_4}
 & \leq c \, (\norma{\nabla\phih}_\infty \, \norma{\nabla\phih}_6 + \norma{\nabla\phistar}_\infty \, \norma{\nabla\phistar}_6)
    \norma{\phih-\phistar}_6^2
  \non
  \\
  & \quad {}
  + c (\norma{\Delta\phih}_6 + \norma{\Delta\phistar}_6) |\norma{\phih-\phistar|}_6^2
  \non
  \\
  & \quad {}
  + c \, (\norma{\nabla\phih}_6 + \norma{\nabla\phistar}_6) \norma{\phih-\phistar}_6 \, \norma{\nabla(\phih-\phistar)}_\infty
  \non
  \\
  & \quad {}
  + c \, \norma{\nabla(\phih-\phistar)}_\infty^2
  + c \, \norma{\phih-\phistar}_4 \, \norma{\Delta(\phih-\phistar)}_4 \,.
  \label{stimaR4}
  \end{align}
  
This ends the preliminary estimate of the most delicate terms,
and we apply it in a while.
As said before, the other terms \juerg{are much easier to estimate}.
\an{Here, we briefly comment on another point. }
As in the proof of Theorem~\ref{Wellposednessl}, we \juerg{still have to}  choose the values of the parameters $M$ and~$N$.
As an example, we consider the last term appearing in \eqref{testprimaF}.
We have
\begin{align}
  & \iO \LLamd \cdot \xxi
  \leq \iO |\mustar| \, \an{|\nabla\psi|\, |\xxi| }
  + \iO |\muh-\mustar| \, |\nabla(\phih-\phistar)| \, |\xxi|
  + \iO |\eta| \, |\nabla\phistar| \, |\xxi|
  \non
  \\
  & \leq \norma\mustar_4 \an{\, \norma{\nabla\psi}\, \norma\xxi_4 }
  + \norma{\muh-\mustar} \, \norma{\nabla(\phih-\phistar)}_4 \, \norma\xxi_4
  + \norma\eta \, \norma{\nabla\phistar}_4 \, \norma\xxi_4
  \non
  \\
  & \leq \frac \alpha{16} \, \normaV\xxi^2
  + c \, \normaV\mustar^2 \, \normaV\psi^2
  + c \, \norma{\muh-\mustar}^2 \, \normaW{\phih-\phistar}^2
  + C^* \, \normaW\phistar^2 \norma\eta^2\,,
  \non
\end{align}
where the special symbol $C^*$ is used in place of the generic~$c$.
By a comparison with the \lhs\ \eqref{lhsF}, 
it is now clear that we need $N$ to be larger than $C^* \, \norma\phistar_{\L\infty W}^2$
(and even much larger, since similar situations appear \juerg{also in other terms}),
\juerg{where we note that} the above norm of $\phistar$ is controlled by the stability inequality~\eqref{stability}.

At the end of the estimates, one recalls all the inequalities that regard the \rhs\ 
corresponding to the \lhs\ \eqref{lhsF} and rearranges.
Then, one integrates over $(0,t)$ with an arbitrary $t\in(0,T]$.
The integral of  \juerg{the} \lhs\ \eqref{lhsF} is essentially given~by
\Beq
  \iot \normaV{\xxi(s)}^2 \, ds
  + \normaV{\psi(t)}^2 
  + \intQt \bigl( |\eta|^2 + |\nabla\Delta\psi|^2 + |\Delta^2\psi|^2 \bigr) \,.
  \label{intlhsF}
\Eeq
Among the integrals that should appear on the corresponding \rhs,
we consider only the one involving $\Delta R_4$.
We have, from \eqref{stimaR4},
\begin{align}
  & \intQt |\Delta R_4|^2
  \leq c \iot \bigl(
    \norma{\nabla\phih(s)}_\infty^2 \, \norma{\nabla\phih(s)}_6^2
    + \norma{\nabla\phistar(s)}_\infty^2 \, \norma{\nabla\phistar(s)}_6^2
  \bigr) \norma{(\phih-\phistar)(s)}_6^4 \, ds
  \non
  \\
  & \quad {}
  + c \iot (\norma{\Delta\phih(s)}_6^2 + \norma{\Delta\phistar(s)}_6^2) \norma{(\phih-\phistar)(s)}_6^4 \, ds
  \non
  \\
  & \quad {}
  + c \iot (\norma{\nabla\phih(s)}_6^2 + \norma{\nabla\phistar(s)}_6^2) \norma{(\phih-\phistar)(s)}_6^2 \, \norma{\nabla(\phih-\phistar)(s)}_\infty^2 \, ds
  \non
  \\
  & \quad {}
  + c \iot \norma{\nabla(\phih-\phistar)(s)}_\infty^4 \, ds
  + c \iot \norma{(\phih-\phistar)(s)}_4^2 \, \norma{\Delta(\phih-\phistar)(s)}_4^2 \, ds 
  \non
  \\[2mm]
  \separa
  & \leq c \, \bigl(
    \norma\phih_{\L2{\Hx3}}^2 \, \norma\phih_{\L\infty W}^2
    + \norma\phistar_{\L2{\Hx3}}^2 \, \norma\phistar_{\L\infty W}^2
  \bigr) \norma{\phih-\phistar}_{\L\infty V}^4
  \non
  \\[2mm]
  & \quad {}
  + c \, \bigl(
    \norma\phih_{\L2{\Hx3}}^2 
    + \norma\phistar_{\L2{\Hx3}}^2 
  \bigr) \norma{\phih-\phistar}_{\L\infty V}^4
  \non
  \\[2mm]
  & \quad {} 
  + c \, \bigl(
    \norma\phih_{\L\infty W}^2 
    + \norma\phistar_{\L\infty W}^2 
  \bigr) \norma{\phih-\phistar}_{\L\infty V}^2 \, \norma{\phih-\phistar}_{\L2{\Hx3}}^2 
  \non
  \\[2mm]
  & \quad {} 
  + c \, \norma{\nabla(\phih-\phistar)}_{\L4\Linfty}^4
  + c \, \norma{\phih-\phistar}_{\L\infty V}^2 \, \norma{\phih-\phistar}_{\L2{\Hx3}}^2 \,,
  \non
\end{align}
and the continuous embedding
$\L\infty V\cap\L2{\Hx3}\emb\L4\Linfty$
can be used to treat the first term on the last line.
Therefore, all the above norms of the difference $\phih-\phistar$
can be estimated \juerg{in terms of} the norm of $\hh$ in $\L2H$, thanks to~\eqref{berlin}.
On the other hand, \juerg{the occurring norms} of $\phih$ and $\phistar$ are controlled by the stability estimate \eqref{stability}.
Hence, we conclude that
\Beq
  \intQt |\Delta R_4|^2 \leq c \, \norma\hh_{\L2H}^4 \,.
  \non
\Eeq

This ends the treatment of the most complicated term that enters the \rhs\ corresponding to the \lhs\ \eqref{intlhsF}.
The others are simpler to \juerg{handle}, \juerg{where} some of them lead to similar inequalities
while the other ones can be dealt with \juerg{using Gronwall's} lemma,
since their coefficients are bounded in~$L^1(0,T)$ thanks to~\eqref{stability}.
We \juerg{eventually can} conclude that
\begin{align}
  & \norma\xxi_{\L2V}
  + \norma\psi_{\L\infty V}
  + \norma\eta_{\L2H}
  \non
  \\
  & \quad {}
  + \norma{\nabla\Delta\psi}_{\L2H}
  + \norma{\Delta^2\psi}_{\L2H}
  \leq c \, \norma\hh_{\L2H}^2 \,.
  \non
\end{align}
By applying elliptic regularity inequalities and a comparison in \eqref{secondaF} and \eqref{terzaF}
(to~recover estimates of $\dt\psi$ and~$\omega$), and recalling the definition \eqref{defX} of~$\calX$,
one obtains~\eqref{frechet}. \juerg{The assertion is thus proved.}
\Edim


\section{Necessary conditions for optimality}
\label{ADJOINT}
\setcounter{equation}{0}

In this section, we give a necessary condition for an admissible control $\ustar$ to be optimal.
The first result  is a simple application of Convex Analysis,
using the convexity of $\Uad$ and the assumptions \Hpcost\ on the cost functional~$\calJ$.
Since both the control-to-state operator $\calS$ and 
the first part $J$ of the cost functional $\calJ$
are \Frechet\ differentiable,
the same is true for the composite map $\Jtilde$ defined in \eqref{defJtilde},
and we can compute its \Frechet\ derivative by the chain rule.
Therefore, we have the following result, which gives a precise form to the variational inequality~\eqref{preNC}.

\Bcor
\label{BadNC}
\juerg{Let} $\ustar\in\Uad$ \juerg{be an optimal control for the control problem~\eqref{control} with associated state
 $\soluzstar:=\calS(\ustar)$}.
Then, there exists \juerg{some} 
$\Lambdastar$ in the subdifferential $\partial G(\ustar)$ such that
\begin{align}
  & b_1 \intQ (\phistar-\phiQ) \phi
  + b_2 \iO (\phistar(T)-\phiO) \phi(T)
  + b_3 \intQ \ustar \cdot (\uu-\ustar)
  \non
  \\
  & \quad {}
  + \intQ \Lambdastar \cdot (\uu-\ustar)
  \geq 0
  \qquad \hbox{for every $\uu\in\Uad$}\,,
  \label{badNC}
\end{align}
where $\phi$ is the second component of the solution $\soluz$ to the 
linearized problem \Pbll\ corresponding to $\hh:=\uu-\ustar$.
\Ecor

This result is not yet satisfactory.
Indeed, it involves the solutions to infinitely many linearized problems, since $\uu$ is arbitrary in~$\Uad$.
To overcome this difficulty, we follow \juerg{the} standard procedure
and introduce the proper adjoint problem associated with $\ustar$ and $\soluzstar$.
\juerg{It} consists in finding a quadruplet $\soluza$ with the regularity
\begin{align}
  & \oomega \in  \pier{\L\infty{\VVz} \cap \L2{\HHx2}},
  \label{regoomega}
  \\
  & p \in \H1{\pier{(\Hx3\cap W)^*}} \cap {\pier{\C0H}} \cap \L2{\pier{\Hx3\cap W}},
  \label{regp}
  \\
  & q \in \pier{\L2V} ,
  \label{regq}
  \\
  & r \in \pier{\L2{V^*}},
  \label{regr}
\end{align}
\Accorpa\Regsoluza regoomega regr
that solves the variational equations
\begin{align}
  & \iO D\oomega : \nabla\zz
  + \iO \lambda(\phistar) \oomega \cdot \zz
  + \iO p \, \nabla\phistar \cdot \zz
  = 0
  \non
  \\
  & \quad \hbox{for every $\zz\in\VVz$ and \aet},
  \label{primaa}
  \\[2mm]
  \separa
  & - \< \dt p , z >
  + \iO \lambda'(\phistar) \vstar \cdot \oomega \, z
  - \iO S'(\phistar) p z
  + \iO f''(\phistar) \wstar q z
  \non
  \\
  & \quad {}
  + \pier{\< r , -\Delta z + f'(\phistar) z>}
  + \iO p \vstar \cdot \nabla z
  - \iO \mustar \oomega \cdot \nabla z
  \non
  \\
  & = \iO \gQ \, z
  \qquad \hbox{for every $z\in \pier{\Hx3\cap W}$ and \aet},
  \label{secondaa}
\end{align}  
\last{where the duality pairing in \eqref{secondaa} is understood between $ (\Hx3\cap W)^*$ and~$\Hx3\cap W$,}
\begin{align}
  & \iO (-\nabla\phistar) \cdot \oomega \, z
  + \iO (-\Delta p) \, z
  \pier{{}=  \iO q z}
  \non
  \\
  & \quad \hbox{for every $z\in H$ and \aet},
  \label{terzaa}
  \\[2mm]
  & \iO \nabla q \cdot \nabla z
  \last{{}+\iO (f'(\phistar)+\nu) q z}
  \pier{{}= \< r , z>}
  \non
  \\
  & \quad \hbox{for every $z\in V$ and \aet},
  \label{quartaa}
\end{align}
and satisfies the final condition
\Beq
  p(T) = \gO \, ,
  \label{cauchya}
\Eeq
\Accorpa\Pbla primaa cauchya
where we have set, for convenience,
\Beq
  \gQ := b_1 (\phistar-\phiQ)
  \aand
  \gO := b_2 \bigl( \phistar(T) - \phiO \bigr) .
  \label{defgQO}
\Eeq
\an{It is worth noting that, in view of \eqref{stability} \pier{and \eqref{hpphiQO}}, it holds that}
\Beq
  \gQ \in \L2H
  \aand
  \gO \in \last{H} \,.
  \label{refgQO}
\Eeq
We have written all the equations in their variational \pier{forms.
However, due to the regularity assumed on the functions that occur,
it is clear that \eqref{terzaa} can be written as a 
partial differential equation on~$Q$.
Namely, we have that 
\Beq
  - \nabla\phistar \cdot \oomega - \Delta p \pier{{}= q{}} 
  \quad \aeQ ,
  \non
\Eeq
so that \eqref{secondaa} may be equivalently replaced by 
\begin{align}
  & - \< \dt p , z >
  + \iO \lambda'(\phistar) \vstar \cdot \oomega \, z
  - \iO S'(\phistar) p z
  + \iO f''(\phistar) \wstar (-\nabla\phistar \cdot \oomega - \Delta p) z
  \non
  \\
& \quad {}
  + \pier{\< r , -\Delta z + f'(\phistar) z>}
  + \iO p \vstar \cdot \nabla z
  - \iO \mustar \oomega \cdot \nabla z = \iO \gQ \, z
  \non
  \\
  & 
  \quad \hbox{for every $z\in \pier{\Hx3\cap W}$ and \aet},
  \label{2piera}
\end{align}
\juerg{and} \eqref{quartaa} becomes
\begin{align}
  & \iO \nabla (- \nabla\phistar \cdot \oomega - \Delta p )\cdot \nabla \zeta
  + \iO (f'(\phistar) +\nu )(- \nabla\phistar \cdot \oomega - \Delta p ) \zeta   \pier{{}= \< r , \zeta>}
  \non
  \\
  & \quad \hbox{for every $\zeta\in V$ and \aet}.
  \label{4piera}
\end{align}
Taking now an arbitrary  $z\in \pier{\Hx3\cap W}$ and letting $\zeta = -\Delta z + f'(\phistar) z $ in \eqref{4piera}, the resulting identity can be used to replace the fifth term on the \lhs\ of~\eqref{2piera}. In this way, we obtain}
\begin{align}
  & - \< \dt p , z >
  + \iO \lambda'(\phistar) \vstar \cdot \oomega \, z
  - \iO S'(\phistar) p z
  + \iO f''(\phistar) \wstar (- \nabla\phistar \cdot \oomega - \Delta p) z
  \non
  \\
  & \quad {}
   \pier{{}+ \iO \nabla (- \nabla\phistar \cdot \oomega) \cdot \nabla(-\Delta z + f'(\phistar) z)
  + \iO \nabla (-\Delta p) \cdot \nabla (-\Delta z +  f'(\phistar) z)}
  \non
  \\
  \separa
  & \quad {}
  - \iO (f'(\phistar)+\nu) (\nabla\phistar \cdot \oomega) (-\Delta z + f'(\phistar) z)
  - \iO (f'(\phistar)+\nu) \Delta p (-\Delta z + f'(\phistar) z)
  \non
  \\
  \separa
  & \quad {}
  + \iO p \vstar \cdot \nabla z
  - \iO \mustar \oomega \cdot \nabla z \pier{{}= \iO \gQ \, z}
  \quad \hbox{for every $z\in V$ and \aet}.
  \label{quintaa}
\end{align}

\an{The first result concerns the well-posedness of} the adjoint problem.

\Bthm
\label{Wellposednessa}
Let $\ustar\in\UR$ be given, and let $\soluzstar:=\calS(\ustar)$ be the corresponding state.
Then, with the notation \eqref{defgQO}, the adjoint problem \Pbla\ has a unique solution satisfying \Regsoluza.
\Ethm

\Bdim
For the sake of brevity, we are forced to perform just formal estimates also in this case.
However, we will make some observations at the end of the proof.
We test \eqref{primaa} by~$\oomega$.
\juerg{From the} coerciveness inequality \eqref{coercive} we then infer~that
\Beq
  \alpha \, \normaV\oomega^2
  \leq - \iO p \, \nabla\phistar \cdot \oomega \,.
  \label{testprimaa}
\Eeq
At the same time, we test \eqref{quintaa} by \pier{$p$} to obtain that\pier{
\begin{align}
  &\pier{{} - \frac 12 \, \ddt \, \norma p^2 + \iO |\nabla \Delta p|^2 }
  \non
  \\
  & 
  = - \iO \lambda'(\phistar) \vstar \cdot \oomega \, \pier{p}
  + \iO S'(\phistar) |\pier{p}|^2
  - \iO f''(\phistar) \wstar (-\nabla\phistar \cdot \oomega - \Delta p) \pier{p}
  \non
  \\
  \separa
  & \quad {}
   \pier{{}- \iO \nabla (- \nabla\phistar \cdot \oomega) \cdot \nabla(-\Delta p + f'(\phistar) p)
  - \iO \nabla (-\Delta p) \cdot \nabla ( f'(\phistar) p)}
  \non
  \\
  \separa
  & \quad {}
  + \iO (f'(\phistar)+\nu) (\nabla\phistar \cdot \oomega) (-\Delta\pier{p} + f'(\phistar) \pier{p})
  \non
  \\
  \separa
  & \quad {}
  + \iO (f'(\phistar)+\nu) \Delta p (-\Delta\pier{p} + f'(\phistar) \pier{p})
  \non
  \\
  \separa
  & \quad {}
  - \iO p \vstar \cdot \nabla \pier{p}
  + \iO \mustar \oomega \cdot \nabla \pier{p}
  + \iO \gQ \, \pier{p} \,.
  \label{testquintaa}
\end{align}
Each term  
on the \rhs\ of \eqref{testquintaa} will be estimated separately.}
In the whole proof, we account for the stability estimate~\eqref{stability} satisfied by the solution~$\calS(\ustar)$.
In particular, $\phistar$ is estimated in $\L\infty W$.
However, \juerg{prior to this}, we make some preliminary observations.
First, we deduce from \eqref{testprimaa}~that
\Beq
  \alpha \, \normaV\oomega^2
  \leq \an{\norma p} \, \norma{\nabla\phistar}\an{_4} \, \norma\oomega_4
  \leq c \, \an{\norma p} \, \normaV\oomega \,,
  \non
\Eeq
whence \juerg{we  immediately conclude that}
\Beq
  \normaV\oomega \leq c \, \an{\norma p} \,.
  \label{stimaoomegaV}
\Eeq
Next, for a given $\gg\in\LLx2$, we consider the problem of finding $\vv\in\VVz$ satisfying the variational equation
\Beq
  \iO (D\vv : \nabla\zz + \vv \cdot \zz)
  = \iO \gg \cdot \zz
  \quad \hbox{for every $\zz\in\VVz$} \,.
  \non
\Eeq
This problem has a unique solution~$\vv$ by the Lax--Milgram theorem.
Moreover, it turns out that $\vv$ belongs to $\HHx2$ and satisfies the estimate
$\norma\vv_{\HHx2}\leq C\,\norma\gg$
with a constant $C$ that depends only on~$\Omega$.
This result can be derived, e.g., from \cite[Lemma 2.49]{E_thesis},
and we apply it to the equation \eqref{primaa} written in the form
\Beq
  \iO (D\oomega : \nabla\zz + \oomega \cdot \zz)
  = \iO \bigl( \an{-} \lambda(\phistar) \oomega \an{-}  p \, \nabla\phistar + \oomega \bigr) \cdot \zz
  \quad \hbox{for every $\zz\in\VVz$} \,,
  \non
\Eeq
for a fixed time.
Owing to \eqref{stimaoomegaV}, we easily see that
\Beq
   \norma{\an{-} \lambda(\phistar) \oomega \an{-}  p \, \nabla\phistar + \oomega}
  \leq c \, \norma\oomega + \norma p_4 \, \norma{\nabla\phistar}_4
  \leq c \, \normaV p\,,
  \non
\Eeq
so that
\Beq
  \norma\oomega_{\HHx2} \leq c \, \normaV p 
  \quad \aet \,.
  \label{stimaoomega}
\Eeq

We can now proceed in \pier{the estimates of the terms on the \rhs\ of \eqref{testquintaa}.}
In the sequel, $\delta$~is a positive parameter.
We owe to the above estimates of $\oomega$ 
and widely account for some of the inequalities \accorpa{sobolev1}{compact2}.
To begin with, we have that
\begin{align}
  \pier{{}-{}} \iO \lambda'(\phistar) \vstar \cdot \oomega \, \pier{p}
  \leq c \, \norma\vstar_4  \, \norma\oomega_4 \norma p  \leq c \, \normaV\vstar\, \norma{\pier{p}}^2 
  \non
\end{align}
\pier{as well as}  
  \begin{align}
  \pier{\iO S'(\phistar) |\pier{p}|^2  \leq c \,\norma{\pier{p}}^2 .}  
  \non
\end{align}
We \juerg{continue with the next term, finding} that
\begin{align}
  & \pier{{}-{}} \iO f''(\phistar) \wstar (-\nabla\phistar \cdot \oomega - \Delta p) \pier{p}
  \non
  \\[1mm]
  & \leq \pier{c \, \norma\wstar_\infty \, \norma{\nabla\phistar}_4 \, \norma\oomega_4 \norma{\pier{p}}
  + c \, \norma\wstar_\infty \, \norma{\Delta p} \, \norma{\pier{p}}}
  \non
  \\[1mm]
  & \leq \pier{c  \norma\wstar_{\Hx2} \norma p^2 + \norma p_{\Hx2}^2  + c \,  \norma\wstar_{\Hx2}^2 \, \norma{\pier{p}}^2 }
  \non
  \\[1mm]
  &\leq  \pier{ \delta \, \norma{\nabla \Delta p}^2 + c_\delta \bigl(1 +   \norma\wstar_{\Hx2}^2 \bigr) \, \norma{\pier{p}}^2   }
 \,.
  \non
\end{align}
The treatment of the next term needs \pier{the little computation
\Beq
  \nabla(\nabla\phistar \cdot \oomega)
  =( D^2\phistar ) \oomega
  + (\nabla\oomega )\nabla \phistar \,,
  \non
\Eeq
where $D^2$ denotes the Hessian operator.
Hence, we have that
\begin{align}
&- \iO \nabla (- \nabla\phistar \cdot \oomega) \cdot \nabla(-\Delta p + f'(\phistar) p)
 \non
  \\
  &=  \iO ( D^2\phistar ) \oomega \cdot \nabla(-\Delta p) +
  \iO  (\nabla\oomega) \nabla \phistar \cdot \nabla(-\Delta p)  
 \non
  \\
  &\quad +   
\iO ((D^2\phistar ) \oomega + (\nabla\oomega) \nabla \phistar) \cdot (f''(\phistar) \last{\nabla \phistar }p + f'(\phistar) \nabla p) 
 \,,
  \non
\end{align}
and we estimate these integrals \juerg{individually. At first, invoking \eqref{stimaoomegaV} and \eqref{stimaoomega},} we infer that 
\begin{align}
  & \iO(D^2\phistar ) \oomega \cdot \nabla(-\Delta p) +
  \iO  (\nabla\oomega ) \nabla \phistar \cdot \nabla(-\Delta p)  
  \non
  \\
  & \leq c \, \norma{\phistar}_{W^{2,4}(\Omega)} \, \norma{\oomega}_4 \, \norma{\nabla \Delta p} 
  + c \, \norma{\nabla \oomega}_{4} \, \norma{\nabla \phistar}_4 \norma{\nabla \Delta p}
  \nonumber
  \\[1mm]
  & \leq \delta \, \norma{\nabla \Delta p}^2 + c_\delta \, \norma{\phistar}_{\Hx3}^2\, \norma{p}^2
   + c_\delta \, \norma{\oomega}_{\HHx2}^2
  \nonumber
  \\
  \separa
  & \leq \delta \, \norma{\nabla \Delta p}^2 + c_\delta \, \norma{\phistar}_{\Hx3}^2\, \norma{p}^2
   + c_\delta \, \norma{p}_{V}^2
   \nonumber
   \\
   & \leq 2\delta \, \norma{\nabla \Delta p}^2 + c_\delta \, \bigl(1+ \norma{\phistar}_{\Hx3}^2\bigr)  \norma{p}^2 \,.
  \non
\end{align}
Similarly, we deduce that
\begin{align*}
&\iO ((D^2\phistar ) \oomega + (\nabla\oomega) \nabla \phistar) \cdot (f''(\phistar)\last{\nabla \phistar }p  + f'(\phistar) \nabla p) 
\\
  & \leq c \, \bigl( \norma{\phistar}_{W^{2,4}(\Omega)} \, \norma{\oomega}_4 + \norma{\nabla \oomega}_{4} \, \norma{\nabla \phistar}_4 \bigr) \bigl( \norma{\nabla \phistar}_4\norma p_4 + \norma{\nabla p} \bigr)
  \\[1mm]
  \separa
    & \leq c \, \bigl( \norma{\phistar}_{\Hx3} \, \norma{p} + \norma{p}_V  \bigr) \norma p_V \leq c \,  \norma{\phistar}_{\Hx3}^2 \, \norma{p}^2 + c \, \norma{p}_V^2
  \\[1mm]
     & \leq \delta \, \norma{\nabla \Delta p}^2 + c_\delta \, \bigl(1+ \norma{\phistar}_{\Hx3}^2 \bigr)  \norma{p}^2 \,.
\end{align*}
Now, about the next term in \eqref{testquintaa}\last{,} we point out that
\begin{align*}
&- \iO \nabla (-\Delta p) \cdot \nabla ( f'(\phistar) p)
=  \iO \nabla (\Delta p) \cdot (f''(\phistar)  \last{\nabla \phistar } p+ f'(\phistar) \nabla p) 
\\
\separa
& \leq c \, \norma{\nabla \Delta p}\,  \bigl( \norma{\nabla \phistar}_4\norma p_4 + \norma{\nabla p} \bigr) 
\\[1mm]
\separa
&\leq \delta \, \norma{\nabla \Delta p}^2 +  c_\delta  \norma p_V^2 
 \leq 2\delta \, \norma{\nabla \Delta p}^2 + c_\delta   \norma{p}^2 \, .
\end{align*}
We continue with the two terms containing the factor $f'(\phistar)+\nu$. Using the same arguments as before, we find
\begin{align}
 &\iO (f'(\phistar)+\nu) (\nabla\phistar \cdot \oomega) (-\Delta\pier{p} + f'(\phistar) \pier{p})
  + \iO (f'(\phistar)+\nu) \Delta p (-\Delta\pier{p} + f'(\phistar) \pier{p})
  \non
  \\
  & \leq c \norma{\nabla\phistar}_4 \, \norma\oomega_4 \, (\norma{\Delta p} + \norma{p})
  + c \, \norma{\Delta p} \, (\norma{\Delta p} + \norma{ p})
  \non
  \\[1mm]
  & \leq c\, \norma{p}_{\Hx2}^2
  \leq  \delta \, \norma{\nabla \Delta p}^2 + c_\delta   \norma{p}^2 \,.
  \non
\end{align}
Next, thanks  to \eqref{stimaoomegaV}\juerg{,} we obtain that
\begin{align}
  & 
  - \iO p \vstar \cdot \nabla \pier{p}
  + \iO \mustar \oomega \cdot \nabla \pier{p}
  \non
  \\
  & \leq c ( \norma p \, \norma\vstar_4 + \norma\mustar \, \norma\oomega_4) \norma{\nabla p}_4
   \non
  \\[2mm]
  & \leq c ( \norma p \, \norma\vstar_V + \norma\mustar \, \norma p) \norma{p}_{\Hx2}
  \non
\\[2mm]
  & \leq  \norma{p}_{\Hx2}^2 + c \bigl(\normaV\vstar^2 + \norma\mustar^2\bigr) \norma p^2 
  \non
  \\[2mm]
  & \leq \delta \, \norma{\nabla \Delta p}^2 + c_\delta  \bigl(1 + \normaV\vstar^2 + \norma\mustar^2\bigr) \norma p^2 \,.
  \non
\end{align}
Finally, we have
\Beq	
  \iO \gQ \, \pier{p} 
  \leq  \frac12 \norma\gQ^2 + \frac12 \norma p^2  \,.
  \non
\Eeq
This ends the preliminary estimates of the terms on the \rhs\
of~\eqref{testquintaa}.}
Thus, we rearrange, choose $\delta>0$ small enough\juerg{,} and integrate with respect to time over $(t,T)$, with an arbitrary $t\in[0,T)$.
\pier{Hence, we obtain
\begin{align}
  &\frac12  \norma{p(t)}^2 
  + \frac12 \intQt |\nabla\Delta p|^2
  \leq  \frac12\last{\norma\gO^2} + \frac12  \intQt |\gQ|^2
  \nonumber
  \\
  &\quad{} 
+ c \int_t^T \Bigl(1 + \normaV{\vstar(s)}^2 + \norma{w^*(s)}_{\Hx2}^2
+ \norma{\phi^*(s)}_{\Hx3}^2 + \norma{\mustar(s)}^2\Bigr) \norma{p(s)}^2 ds  \,,
  \non
\end{align}
where now $Q_t:=\Omega\times(t,T)$. Since the function 
$$
s\mapsto 1 + \normaV{\vstar(s)}^2 + \norma{w^*(s)}_{\Hx2}^2
+ \norma{\phi^*(s)}_{\Hx3}^2 + \norma{\mustar(s)}^2
$$
is bounded in $L^1(0, T)$ by virtue of~\eqref{stability}, 
the (backward) Gronwall lemma allows us to infer that
\begin{align}
  & \norma p_{\L\infty H}
  + \norma{\nabla\Delta p}_{\L2 H}
  \leq c \bigl( \norma\gQ_{\L2H} + \last{\norma\gO} \bigr) \,.
  \non
\end{align}
This estimate is improved using the estimate~\eqref{elliptic1.5}
from elliptic regularity theory and, consequently, the inequalities~\eqref{stimaoomegaV} and~\eqref{stimaoomega}, in order to conclude that
\begin{align}
  & 
  \last{\norma p_{\L\infty H\cap \L2{\Hx3\cap W}} + \norma{\oomega}_{\L\infty {\VV_0}\cap \L2{\HHx2} }}
  \nonumber
  \\
  &
  \leq c \bigl( \norma\gQ_{\L2H} + \last{\norma \gO} \bigr) \,.
  \label{fin-a}
\end{align}
Then, additional} estimates for $\dt p$, $q$ and $r$ are easily obtained via comparison in
\eqref{secondaa}, \eqref{terzaa} and~\eqref{quartaa}, respectively.
We thus deduce~that
\begin{align}
  &\pier{\norma p_{\H1{\pier{(\Hx3\cap W)^*}}}
  + \norma q_{\L2V}
  + \norma r_{\L2{V^*}}}
  \non
  \\
  & \leq c \bigl( \norma\gQ_{\L2H} + \last{\norma \gO} \bigr) \,.
  \label{finea}
\end{align}
We conclude with some observation for the reader's convenience.
First, the procedure that led us to \pier{the estimates~\eqref{fin-a} and} \eqref{finea} is \juerg{rigorous} if performed on a solution satisfying \Regsoluza.
In the particular case when $\gQ=0$ and $\gO=0$, this estimate implies that the solution vanishes.
This proves the uniqueness of the solution by linearity.
As for existence, the formal estimates we have derived can be performed 
on an approximating problem based on a Faedo--Galerkin scheme by assuming that
$\gQ$ and $\gO$ are given by~\eqref{defgQO}.
\an{Then, by standard compactness arguments, one can construct a solution possessing the asserted regularity.}
\Edim

Our \juerg{next result yields the following first-order necessary optimality condition}.

\Bthm
\label{GoodNC}
Assume that $\ustar\in\Uad$ is an optimal control for the control problem~\eqref{control},
and let $\soluzstar:=\calS(\ustar)$ be the corresponding state.
Then there exists some $\Lambdastar$ in the subdifferential $\partial G(\ustar)$ such that
\Beq
  \intQ \bigl( b_3 \ustar + \oomega \bigr) \cdot (\uu-\ustar)
  + \intQ \Lambdastar \cdot (\uu-\ustar)
  \geq 0
  \qquad \hbox{for every $\uu\in\Uad$}\,,
  \label{goodNC}
\Eeq
where $\oomega$ is the first component of the solution $\soluza$ to the adjoint problem \Pbla.
\Ethm

\Bdim
We fix $\uu\in\Uad$ and consider the linearized problem \Pbll\ corresponding to $\hh:=\uu-\ustar$.
\last{Let us emphasize that \eqref{quartal} is equivalent to \eqref{4pierl}, whose terms are all in $V$, \aet.} 
We test \pier{\eqref{primal} by $\oomega$, \eqref{secondal} by $p$, \eqref{terzal} by $q$, and take the duality pairing between $r$ and \eqref{4pierl}, respectively. Then, we}
integrate over $(0,T)$ with respect to time 
and take the sum of the resulting equalities.
\juerg{It then follows the identity}
\begin{align}
  & \intQ D\vv : \nabla\oomega
  + \intQ \lambda(\phistar) \vv \cdot \oomega
  + \intQ \lambda'(\phistar) \phi \vstar \cdot \oomega
  - \intQ \bigl(
    \mu \nabla\phistar
    + \mustar \nabla\phi
    + \uu - \ustar
  \bigr) \cdot \oomega
  \non
  \\
  & \quad {}
  + \ioT \< \dt\phi(t) , p(t) > \, dt
  + \intQ \vv \cdot \nabla\phistar \, p
  + \intQ \vstar \cdot \nabla\phi \, p
  + \intQ \mu (-\Delta p)
  - \intQ S'(\phistar) \phi p
  \non
  \\
  & \quad {}
  + \intQ \nabla w \cdot \nabla q
  + \intQ f''(\phistar) \phi \wstar q
  + \intQ f'(\phistar) w q 
  + \nu \intQ w q
  - \intQ \mu q
  \non
  \\
  & \quad {}
  + \pier{\int_0^T \<r(t), (-\Delta\phi +  f'(\phistar) \phi  - w) (t) >\,  dt}= 0 \,.
  \label{testlin}
\end{align}
At the same time, we consider the adjoint problem \Pbla\ and test its equations by $\vv$, $\phi$, $\mu$ and~$w$,
integrate over $(0,T)$ and sum up.
We obtain
\begin{align}
  & \intQ D\oomega : \nabla\vv
  + \intQ \lambda(\phistar) \oomega \cdot \vv
  + \intQ p \, \nabla\phistar \cdot \vv
  \non
  \\
  & \quad {}
  - \ioT \< \dt p(t) , \phi(t) > \, dt
  + \intQ \lambda'(\phistar) \vstar \cdot \oomega \, \phi
  - \intQ S'(\phistar) p \phi
  + \intQ f''(\phistar) \wstar q \phi
  \non
  \\
  \separa
  & \quad {}
  + \pier{\int_0^T \<r(t), (-\Delta\phi +  f'(\phistar) \phi) (t) >\,  dt}
  + \iO p \vstar \cdot \nabla\phi
  - \intQ \mustar \oomega \cdot \nabla\phi
  - \intQ \gQ \, \phi
  \non
  \\
  \separa
  & \quad {}
  + \intQ (-\nabla\phistar) \cdot \oomega \, \mu
  + \intQ (-\Delta p) \, \mu
  - \intQ q \mu
  \non
  \\
  & \quad {}
  + \intQ \nabla q \cdot \nabla w
  + \intQ f'(\phistar) q w 
  + \nu \intQ q w
  - \pier{\int_0^T \<r(t), w (t) >\,  dt}
  = 0 \,.
  \label{testadj}
\end{align}
At this point, we take the difference between \eqref{testlin} and \eqref{testadj}.
Several terms cancel each other, and it remains the following identity:
\Beq
  - \intQ (\uu - \ustar) \cdot \oomega
  + \ioT \< \dt\phi(t) , p(t) > \, dt
  + \ioT \< \dt p(t) , \phi(t) > \, dt
  + \intQ \gQ \phi
  = 0 \,.
  \non
\Eeq
Since both $\phi$ and $p$ belong to \pier{$\H1{\pier{(\Hx3\cap W)^*}}\cap\L2{\Hx3\cap W}$},
the sum of the two integrals in the middle of the above line
is the time integral of the time derivative of~$\iO\phi\,p$.
Thus, on account of the initial condition \eqref{cauchyl}, the final condition \eqref{cauchya}
and the definitions \eqref{defgQO}, we deduce that
\Beq
  - \intQ (\uu - \ustar) \cdot \oomega
  + \iO  b_2 \bigl( \phistar(T) - \phiO \bigr) \, \pier{\phi(T)} 
  + \intQ b_1 (\phistar-\phiQ) \, \phi
  = 0 \,.
  \non
\Eeq
By combining this with the inequality \eqref{badNC},
we obtain \eqref{goodNC}.
\Edim

\vspace{2mm}
\Brem
\label{Oss1}
The variational inequality \eqref{goodNC} is equivalent to a system of three 
variational inequalities that have to be satisfied simultaneously. Indeed, if 
$\Lambdastar=(\Lambda_1^*,\Lambda_2^*,\Lambda_3^*)$ is given, then 
it is easily seen that
\eqref{goodNC} holds true if and only if the components of $\uu^*
=(u_1^*,u_2^*,u_3^*)$ satisfy
for $i=1,2,3$ the variational inequalities
\begin{align}
\label{goodNCi}
\int_Q (b_3 u^*_i+\omega_i+\Lambda_i^*)(u-u_i^*) \,\ge\,0  \quad\mbox{for every }\,
u\in U_{ad}^i\,,
\end{align}
where $\oomega=(\omega_1,\omega_2,\omega_3)$ and
\begin{equation}
\label{Uadi}
U_{ad}^i :=\{u\in L^\infty(Q):  \underline u_i \le u\le \overline u_i \,\mbox{ a.e. in }\,Q\},
\quad i=1,2,3.
\end{equation} 
A standard argument then shows that the following pointwise projection conditions
are valid:
\begin{equation}
\label{pro1}
u_i^*\,=\,\max \left\{\underline u_i, \,\min \left\{\overline u_i, - b_3^{-1}(\omega_i
+\Lambda_i^* \last{)}\right\} \right\} \quad\mbox{a.e. in \,$Q$, for }\,i=1,2,3.
\end{equation}
\Erem
\Brem
\label{Oss2}
It is worth noting that the argumentation used to show the validity 
 of Corollary~\ref{BadNC}, Theorem~\ref{GoodNC}, and Remark~\ref{Oss1} 
does not only work for (globally) optimal controls, but also for controls that are only
 locally optimal in the sense of $\LL^\infty(Q)$. The results shown above 
 therefore hold true correspondingly in the locally optimal case.
  We recall in this connection that a
  control  $\uu^*\in \Uad$ is \last{termed} locally optimal in the sense of $\LL^p(Q)$ for 
  some $p\in[1,+\infty]$ if there is some $\gamma>0$ such that
  \Beq
  \label{locopt}
  {\cal J}(\phi^*,\ustar)\,\le\,{\cal J}(\phi,\uu) \quad\mbox{for all $\uu\in \Uad$ with }\,
  \|\uu-\ustar\|_{\LL^p(Q)}\le\gamma\,,
  \Eeq
  where $\phi$ and $\phi^*$ denote the second components of ${\cal S}(\uu)$
  and ${\cal S}(\ustar)$, respectively. 
  Observe that every control, which is locally optimal in the sense of $\LL^p(Q)$ for
  some $p\in [1,+\infty)$, is also locally optimal in the sense of $\LL^\infty(Q)$.  
\Erem

\section{Remarks on sparsity}
\label{SPARSITY}
\setcounter{equation}{0}

In this section, we discuss the aspect of sparsity in the optimal control problem under
investigation. We again assume that the general assumption \eqref{allassumptions} is fulfilled
and that $\ustar\in\Uad$ is an optimal control with associated state $(\vv^*,\phi^*,\mu^*,w^*)
={\cal S}(\ustar)$ and adjoint state $(\oomega,p,q,r)$. We remark at this place that all of the following results would remain valid correspondingly if $\ustar$ were only locally optimal in the sense of $\LL^\infty(Q)$. 

The convex functional $G$ in the objective functional is \last{responsible} for the occurrence of
sparsity, i.e., the possibility that optimal controls may vanish in some subregion of $Q$
having a positive measure. The special sparsity properties are determined by the form
of the subdifferential $\partial G$  via the variational inequality \eqref{goodNC}. 
There are several concepts of sparsity, which are each induced by a specific functional $G$.
We confine ourselves to {\em full sparsity}, which is induced by the $L^1(Q)$ norm
\Beq
\label{defj}
j(u):=\|u\|_{L^1(Q)}=\int_Q |u|\,, \quad\mbox{for $u\in L^2(Q)$ }, 
\Eeq
which is nonnegative, convex and continuous (and thus sequentially lower semicontinuous)
on the space $L^2(Q)$. It is well known (see, e.g., \cite{Ioffe}) that 
its subdifferential $\partial j(u)$ is for 
every $u\in L^2(Q)$ given by the set of all $\lambda\in L^2(Q)$ that satisfy, for almost
every $(x,t)\in Q$,
\begin{equation}\label{partialj}
\pier{\lambda(x,t)\in
\begin{cases}
\{+1\},   & \text{if } u(x,t)>0,\\[4pt]
[-1,+1],  & \text{if } u(x,t)=0,\\[4pt]
\{-1\},   & \text{if } u(x,t)<0.
\end{cases}}
\end{equation}
In the following, we consider a sparsity term which is slightly more general than that given
in \eqref{Gsparse}, namely
\begin{align}
\label{Gneu}
&G(\uu)=\int_Q\bigl(\kappa_1|u_1|+\kappa_2|u_2|+\kappa_3|u_3|\bigr)
=\sum_{i=1}^3 \kappa_i j(I_i(\uu)) \non\\
&\quad{} \mbox{for }\, \uu=(u_1,u_2,u_3)\in \LL^2(Q)\,,
\end{align}
with given sparsity parameters $\kappa_i>0$, $i=1,2,3$, and the linear and continuous projection 
operators $I_i: \LL^2(Q)\to L^2(Q), \uu=(u_1,u_2,u_3)\mapsto u_i, \,\,\,i=1,2,3.$ Using the
well-known rules for subdifferentials (see, e.g., \cite[Sect.~4.2.2, Thms.~1~and~2]{Ioffe}),
and denoting by $I_i^*$ the dual operator of $I_i$, $i=1,2,3$,
we then conclude that the subdifferential of $G$ is given~by
\begin{align}
\partial G(\uu)&\,=\,\sum_{i=1}^3 \kappa_i \,I_i^*\,\partial j(I_i(\uu))\non\\
&\,=\,\{(\kappa_1\lambda_1,\kappa_2 \lambda_2,\kappa_3\lambda_3)\in\LL^2(Q):
\lambda_i\in \partial j(u_i) , \,\,\,i=1,2,3\}.
\end{align}
Therefore, the components of the multiplier $\Lambdastar=(\Lambda_1^*,\Lambda_2^*,
\Lambda_3^*)\in\partial G(\ustar)$ occurring in \eqref{goodNC} are\last{,} for $i=1,2,3$\last{,} of the form
$\,\Lambda_i^*=\kappa_i\lambda_i^*$ with some $\lambda_i^*\in\partial j(u_i^*)$.
From \eqref{goodNCi} and \eqref{pro1} in Remark~\ref{Oss1} we then infer that for $i=1,2,3$ the variational inequality
\Beq
\int_Q \bigl(b_3 u_i^* + \omega_i + \kappa_i \lambda_i^*\bigr)(u-u_i^*) \,\ge\,0
\quad\mbox{for every }\,u\in U_{ad}^i
\Eeq
and the projection condition
\Beq
\label{pro2}
u_i^*\,=\,\max 
\left \{ \underline u_i, \,\min \left \{\overline u_i, -b_3^{-1}(\omega_i
+\kappa_i \lambda_i^*) 
\right \} \right \} \quad\mbox{a.e. in  }\,Q\,,
\Eeq
are satisfied, with suitable $\,\lambda_i^*\in\partial j(u_i^*)$.
We then \pier{arrive at the following sparsity result, which is in line with analogous results established} 
in~\cite{CoSpTr,SprTro,SpTr}\last{.}
\Bthm
\label{Sparse}
Suppose that the general assumptions for our optimal control problem are satisfied, and 
assume that the thresholds satisfy the following condition:
\Beq
\label{LagoMaggiore}
\underline u_i \mbox{\, and }\, \overline u_i \,\mbox{ are constants such that }\,
\underline u_i < 0 <\overline u_i\,\mbox{ for }\,i=1,2,3. 
\Eeq
Moreover, let $\ustar\in\Uad$ be an optimal control with associated state 
$(\vv^*,\phi^*,\mu^*,w^*)={\cal S}(\ustar)$ and adjoint state $(\oomega,p,q,r)$.
Then,  for every $i\in\{1,2,3\}$, it holds  that for a.e. $(x,t)\in Q$ the following equivalence  
is valid:
\Beq
\label{LagodiComo}
u_i^*(x,t)=0   \quad\Longleftrightarrow\quad |\omega_i(x,t)|\le \kappa_i\,.
\Eeq
\Ethm
\Bdim
Let $i\in\{1,2,3\}$ be fixed. For almost every $(x,t)\in Q$ it holds: if $u_i^*(x,t)=0$, then,
by virtue of \eqref{pro2}, $0=-b_3^{-1}(\omega_i(x,t)+\kappa_i \lambda_i^*(x,t))$, whence 
we immediately obtain  that 
\,$|\omega_i(x,t)|=\kappa_i |\lambda_i^*(x,t)|$. Since $\lambda_i^*\in\partial j(u_i^*)$, 
it then follows from \eqref{partialj} that \,$|\omega_i(x,t)|\le \kappa_i$.

Conversely, we have for almost every $(x,t)\in Q$ \pier{the following implication: if}  $|\omega_i(x,t)|\le\kappa_i$\, and 
\,$u_i^*(x,t)>0$, then, by \eqref{partialj}, \pier{it follows that} $\,\lambda_i^*(x,t)=+1\,$ and thus, again by \eqref{pro2},   
$0<-b_3^{-1}(\omega_i(x,t)+\kappa_i)$. But then 
$$ \omega_i(x,t)+\kappa_i<0\, \hbox{ and therefore 
\,$|\omega_i(x,t)|=-\omega_i(x,t)>\kappa_i$,}$$ 
\last{leading to} a contradiction. Similar reasoning yields
a contradiction also if $\,u_i^*(x,t)<0\,$ is assumed. In conclusion, we must have for almost
every $(x,t)\in Q$: if \,$|\omega_i(x,t)|\le\kappa_i$\, then $\,u_i^*(x,t)=0$. This
 concludes the proof of the assertion.
\Edim
\Brem
\label{Oss3}
Observe that the adjoint variable $\oomega$ appearing in the sparsity condition 
\eqref{LagodiComo} depends on the special optimal control $\ustar$. It is therefore
natural to rise the question whether there exists a global sparsity parameter $\kappa^*>0$ 
such that all optimal controls vanish  a.e.~in $Q$ whenever the parameters $\kappa_i$ 
exceed the value $\kappa^*$. \juerg{Apparently, this requires to establish a global  
$\LL^\infty(Q)$\last{-}bound for the adjoint variable $\oomega$. 
However,  recalling the global estimate (2.33) established in Theorem 2.1, 
we can conclude from a closer inspection of the estimates performed in the proof of Theorem 5.2
that there exists a global constant $\widehat C>0$ such that
$$
\|\oomega\|_{L^\infty(0,T;\VVz)} \le  \widehat C\,,
$$ 
whenever $\oomega$ is the first component of the solution to the adjoint system associated with an arbitrary optimal control  $\ustar$. Hence, the existence of a suitable global sparsity parameter $\kappa^*$ can be guaranteed at least in the spatially one-dimensional case.}
\Erem


\section*{\pier{Acknowledgments}}

\noindent
\pier{PC and AS benefited from stimulating discussions related to this work during their stay at the Erwin Schr\"odinger International Institute for Mathematics and Physics (ESI) in Vienna, and gratefully acknowledge the warm hospitality of the Institute. PC and AS also acknowledge their affiliation with the GNAMPA (Gruppo Nazionale per l'Analisi Matematica, la Probabilit\`a e le loro Applicazioni) of INdAM (Isti\-tuto Nazionale di Alta Matematica) and the financial support of the GNAMPA project (CUP~E5324001950001) of INdAM. PC further acknowledges partial support from the Next Generation EU Project No.~P2022Z7ZAJ (A unitary mathematical framework for modelling muscular dystrophies).}
\an{AS also acknowledges partial support from the ``MUR Grant -- Dipartimento di Eccellenza'' 2023-2027 and from the Alexander von Humboldt Foundation.}

\section*{Declarations}

\subsection*{Conflicts of Interest}

The authors declare that they have no conflicts of interest and that no competing interests exist regarding the publication of this manuscript.


\End{document}

\section*{Appendix}
\an{
\begin{proposition}\label{APP:ineq}
It holds that, for every $x,y\in\erre^3$,
\begin{align}
	\bigl||x|^2x-|y|^2y\bigr|\leq \tfrac 32 |x-y|(|x|^2+|y|^2).
 \label{inequality}
\end{align}
\end{proposition}
\begin{proof}
First, observe the identity
\begin{align*}
|x|^2x - |y|^2y 
&= |x|^2(x-y) + (|x|^2 - |y|^2)y .
\end{align*}
Taking absolute values to both sides yields
\begin{align*}
\big| |x|^2x - |y|^2y \big|
&\le |x|^2|x-y| + \big||x|^2 - |y|^2\big|\,|y| .
\end{align*}
Since
\begin{align*}
\big||x|^2 - |y|^2\big|
&= \big||x|-|y|\big|\,\big(|x|+|y|\big)
\le |x-y|(|x|+|y|) ,
\end{align*}
we obtain
\begin{align*}
\big| |x|^2x - |y|^2y \big|
&\le |x-y|\big(|x|^2 + |x||y| + |y|^2\big) .
\end{align*}
We now look for a constant $C^*>0$ such that
\begin{align*}
|x-y|\big(|x|^2 + |x||y| + |y|^2\big)
&\le C^*|x-y|\big(|x|^2 + |y|^2\big),
\end{align*}
that is,
\begin{align*}
\frac{|x|^2 + |x||y| + |y|^2}{|x|^2 + |y|^2}
&\le C^* .
\end{align*}
Let $a := |x|$ and $b := |y|$ and oserve that in the case $a=0$ the above inequality is trivially fulfilled. 
Hence, assume 
$a \not =0$ and set $z := \frac{b}{a}$. Consider the nonnegative function
\begin{align*}	
	g(z):= \frac{1 + z + z^2}{1 + z^2}=\frac{a^2 + ab + b^2}{a^2 + b^2}.
\end{align*} 
We now determine the maximum of $g$. A direct computation gives $
g'(z) = \frac{1 - z^2}{(1+z^2)^2},$
so $g$ attains its maximum at $z = 1$, and
$g(1) = \frac{3}{2}$.
Therefore, the optimal constant is $C^* = \frac{3}{2}$.
\end{proof}
}